\documentclass[twoside, a4paper]{amsart}
\usepackage{hyperref} 
\usepackage{amsmath, amssymb, enumerate}
\usepackage[table]{xcolor}
\usepackage{amsaddr}
\usepackage{multirow}
\usepackage{footnote}
\usepackage{booktabs, caption}
\usepackage{pgfplots}\pgfplotsset{compat=1.15}
\usepackage{graphicx, subcaption}
\usepackage{amsthm}
\newtheorem{theorem}{Theorem}[section]
\newtheorem{lemma}[theorem]{Lemma}

\newtheorem{prop}[theorem]{Proposition}
\theoremstyle{definition}
\newtheorem{remark}[theorem]{Remark}
\newtheorem{example}[theorem]{Example}
\newtheorem{definition}[theorem]{Definition}

\usepackage{algorithm,algpseudocode}
\newcommand{\set}[2]{\ensuremath{\{#1 \,\colon #2\}}}
\newcommand{\evalat}[1]{\bigr\rvert_{#1}}
\providecommand{\abs}[1]{\left\lvert{#1}\right\rvert}
\newsavebox{\ldecpartsymb}\sbox{\ldecpartsymb}{$\{\mkern-6mu\{$}
\newsavebox{\rdecpartsymb}\sbox{\rdecpartsymb}{$\}\mkern-6mu\}$}
\makeatletter
\newcommand{\decpart}[1]{\setbox\z@\hbox{$#1$}%
     \ifdim\ht\z@<\ht\ldecpartsymb\usebox\ldecpartsymb\usebox\z@\usebox\rdecpartsymb%
     \else\resizebox{\wd\ldecpartsymb}{1.15\ht\z@}{\usebox\ldecpartsymb}%
          \usebox\z@\resizebox{\wd\rdecpartsymb}{1.15\ht\z@}{\usebox\rdecpartsymb}%
     \fi%
}
\def\@map#1#2[#3]{\mbox{$#1 \colon #2 \longrightarrow #3$}}
\def\map#1#2{\@ifnextchar [{\@map{#1}{#2}}{\@map{#1}{#2}[#2]}}
\makeatother
\newcommand{\includeepssubfigure}[2]{\begin{subfigure}[t]{0.32\textwidth}\centering
  \includegraphics[width=0.7\linewidth, angle=-90]{Images/#1-eps-converted-to.pdf}
   \captionsetup{width=0.8\linewidth}\caption{#2}
\end{subfigure}}
\newcommand{\includepngsubfigure}[2]{\begin{subfigure}[t]{0.32\textwidth}\centering
  \includegraphics[width=0.99\linewidth]{Images/#1.png}
   \captionsetup{width=0.8\linewidth}\caption{#2}
\end{subfigure}}
\DeclareMathOperator{\rot}{Rot}
\DeclareMathOperator{\Orb}{Orb}
\DeclareMathOperator{\diam}{diam}
\newcommand{\N}{\ensuremath{\mathbb{N}}}
\newcommand{\Z}{\ensuremath{\mathbb{Z}}}
\newcommand{\Q}{\ensuremath{\mathbb{Q}}}
\newcommand{\R}{\ensuremath{\mathbb{R}}}
\newcommand{\SI}{\ensuremath{\mathbb{S}^1}}
\newcommand{\dol}[1][1]{\ensuremath{\mathcal{L}_{#1}}}
\title{An Algorithm to compute Rotation Numbers in the circle}
\thanks{Supported by the Spain's ``Agencial Estatal de Investigaci\'on'' (AEI) grants
MTM2017-86795-C3-1-P and MDM-2014-0445 within the ``Mar\'{\i}a de Maeztu'' Program.}
\author{Llu\'is Alsed\`a}
\address{Departament de Matem\`atiques,
  Edifici Cc,
  Universitat Aut\`onoma de Barcelona,
  08193 Bellaterra (Barcelona),
  Spain
}
\address{Centre de Recerca Matem\`atica,
  Campus de Bellaterra,
  Edifici Cc,
  Universitat Aut\`onoma de Barcelona,
  08193 Bellaterra (Barcelona),
  Spain
}
\email{alseda@mat.uab.cat, llalseda@crm.cat}
\author{Salvador Borr\'os-Cullell}
\address{Departament de Matem\`atiques,
  Edifici Cc,
  Universitat Aut\`onoma de Barcelona,
  08193 Bellaterra (Barcelona),
  Spain
}
\email{sborros@mat.uab.cat}
\makeatletter
\@namedef{subjclassname@2020}{\textup{2020} Mathematics Subject Classification}
\makeatother
\subjclass[2020]{Primary 37E10, Secondary 37E45}
\date{December 6, 2020}
\keywords{Rotation number, circle maps, nondecreasing degree one lifting, algorithm}
\begin{document}
\begin{abstract}
In this article we present an efficient algorithm to compute rotation
intervals of circle maps of degree one. It is based on the
computation of the rotation number of a monotone circle map of degree
one with a constant section. The main strength of this algorithm is
that it computes \emph{exactly} the rotation interval
of a natural subclass of the continuous non-invertible degree one circle maps.

We also compare our algorithm with other existing ones by plotting the
Devil's Staircase of a one-parameter family of maps and the Arnold Tongues
and rotation intervals of some special non-differentiable families,
most of which were out of the reach of the existing algorithms that were
centred around differentiable maps.
\end{abstract}
\maketitle
\section{Introduction}\label{sec:Introduction}
The rotation interval plays an important role in combinatorial dynamics.
For example Misiurewicz's Theorem~\cite{misiuthm} links the
set of periods of a continuous lifting $F$ of degree one to the set
$
M := \set{n\in\N}{\tfrac{k}{n} \in \rot(F)\text{ for some integer }k},
$
where $\rot(F)$ denotes the rotation interval of $F.$ Moreover, it is
natural to compute lower bounds of the topological entropy depending
on the rotation interval \cite{entcirc}. In any case, the knowledge
of the rotation interval of circle maps of degree one is of
theoretical importance.

The rotation number was introduced by H.~Poincaré to study the
movement of celestial bodies \cite{poincare}, and since then has been
found to model a wide variety of physical and sociological processes.
The application to voting theory \cite{vots2, vots1} is specially
surprising in this context.

The computation of the rotation number for invertible maps of
degree 1 from $\SI$ onto itself is well studied, and many very
efficient algorithms exist for its computation
\cite{PMIHES_1979__49__5_0,Pavani, Seara-Villanueva, VanVeldhuizen}.
However, there is a lack of an efficient algorithm for the non-invertible
and non-differentiable case.

In this article, we propose a method that allows us to compute
the rotation interval for the non-invertible case.
Our algorithm is based on the fact that we can compute \emph{exactly}
the rotation number of a natural subclass of the of the class of
continuous non-decreasing degree one circle maps that have a constant
section and a rational rotation number.
From this algorithm we get an efficient way to compute
\emph{exactly} the rotation interval of a natural subclass of the
continuous non-invertible degree one circle maps by using the so
called \emph{upper} and \emph{lower} maps,
which, when different, always have a constant section.

To check the efficiency of our algorithm will use it to compute some
classical results such as a Devil's Staircase.
When doing so, we will compare the efficiency of our algorithm with the
performance of some other algorithms that have been traditionally used
under the hypothesis of non-invertibility. 
On the other hand, we will also compute the rotation interval and
Arnold tongues for a variety of maps, in the same comparing spirit.
These maps include the Standard Map and variants of it but have issues
either with the differentiability, or even with the continuity.
Of course these variants are not well suited for algorithms that
strongly use differentiability.

The paper is organised as follows.
In Section~\ref{sec:prelims} the theoretical background will be set.
In Section~\ref{sec:methods} the algorithm will be presented, and
in Section~\ref{sec:examples} we will provide the mentioned examples
of the use of the algorithm.
Finally in Section~\ref{sec:Conclusions} we will discuss the advantages
and disadvantages of the proposed algorithm.
\section{A short Survey on Rotation Theory and the Computation of Rotation Numbers}\label{sec:prelims}

We will start by recalling some results from the rotation theory for
circle maps. To do this we will follow \cite{ALM-Book}.

The \emph{floor function}
(i.e. the function that returns the greatest integer less than or equal to the variable)
will be denoted as $\lfloor{}\thinspace\cdot\thinspace\rfloor.$
Also the \emph{decimal part} of a real number $x \in \R,$ defined as
$x - \lfloor x \rfloor \in [0,1)$ will be denoted by $\decpart{x}.$

In what follows $\SI$ denotes the circle, which is defined as the set of
all complex numbers of modulus one.
Let $\map{e}{\R}[\SI]$ be the natural projection from $\R$ to $\SI,$
which is defined by  $e(x) := \exp(2\pi{}ix).$

Let $\map{f}{\SI}$ be continuous map. A continuous map $\map{F}{\R}$ is a
\emph{lifting of $f$} if and only if
$
   e(F(x)) = f(e(x)) \text{ for every $x\in{}\R$.}
$
Note that the lifting of a circle map is not unique, and that any two liftings
$F$ and $F'$ of the same continuous map $\map{f}{\SI}$ verify
$F = F' + k$ for some $k\in{}\Z.$

For every continuous map $\map{f}{\SI}$ there exists an integer $d$ such that
\[
    F(x+1) = F(x) + d
\]
for every lifting $F$ of $f$ and every $x\in{}\R$
(that is, the number $d$ is independent of the choice of the lifting
and the point $x\in \R$).
We shall call this number $d$ \emph{the degree of $f$}.
The degree of a map roughly corresponds to the number of times that the whole image of the
map $f$ covers homotopically $\SI.$

In this paper we are interested studying maps of degree~1, since the
rotation theory is well defined for the liftings of these maps.

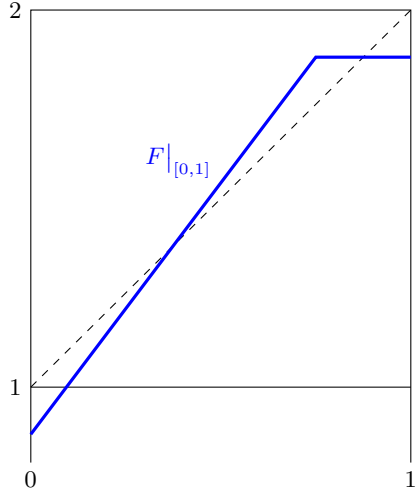
\begin{figure}
\centering
\begin{tikzpicture}[scale=5]
\draw (0,-0.2) -- (0, 1) -- (1, 1) -- (1, -0.2);
\draw (0, 0) -- (1, 0);
\draw[dashed] (0, 0) -- (1, 1);

\draw[blue, very thick] (0, -0.125) -- (0.75, 0.875) -- (1, 0.875);

\node[below] at (0, -0.2) {\small $0$};
\node[below] at (1, -0.2) {\small $1$};
\node[left] at (0,0) {\small $1$};
\node[left] at (0,1) {\small $2$};
\node[left, blue] at (0.5, 0.6) {\small $F\evalat{[0,1]}$};
\end{tikzpicture}
\caption{An example of a map from $\dol$ which can be considered
as a toy model for the elements of that class.
The picture shows $F\evalat{[0,1]},$ and $F$ is globally defined as
$F(x) = F\evalat{[0,1]}\bigl(\decpart{x}\bigr) + \lfloor x \rfloor.$}\label{mupicture}
\end{figure}

We will denote the set of all liftings of maps of degree 1 by $\dol.$
Observe that to define a map from $\dol$ it is enough to define
$F\evalat{[0,1]}$ (see Figure~\ref{mupicture}) because $F$ can be
globally defined as
$F(x) = F\evalat{[0,1]}\bigl(\decpart{x}\bigr) + \lfloor x \rfloor$
for every $x \in \R.$

\begin{remark}\label{iterats-de-grau-1}
It is easy to see that, for every $F\in \dol$,
$F^n(x+k) = F^n(x)+k$ for every
$n\in\N,$ $x\in\R$ and $k\in{}\Z.$
Consequently, $F^n \in \dol$ for every $n\in\N.$
\end{remark}

\begin{definition}\label{defrotnum}
Let $F\in{}\dol$, and let $x\in{}\R.$
We define the \emph{rotation number of $x$} as
\[
\rho_{_F}(x) := \limsup_{n\to\infty}\frac{F^n(x)-x}{n}.
\]
Observe (Remark~\ref{iterats-de-grau-1}) that,
$\rho_{_F}(x) = \rho_{_F}(x+k)$ for every $k\in{}\Z.$
The \emph{rotation set of $F$} is defined as:
\[
   \rot(F) = \set{\rho_{_F}(x)}{x\in{}\R}
           = \set{\rho_{_F}(x)}{x\in[0,1]}.
\]
Ito~\cite{ito}, proved that the rotation set is a closed interval
of the real line. So, henceforth the set $\rot(F)$ will be called
the \emph{rotation interval of $F$}.
\end{definition}

\begin{prop}\label{nondecreasing}
Let $F\in{}\dol$ be non-decreasing.
Then, for every $x\in{}\R$ the limit
\[
   \lim_{n\to\infty}\frac{F^n(x)-x}{n}
\]
exists and is independent of $x.$
\end{prop}

For a non-decreasing map $F\in{}\dol,$ the number
$\rho_{_F}(x) = \lim_{n\to\infty}\frac{F^n(x)-x}{n}$
will be called the \emph{rotation number of $F$},
and will be denoted by $\rho_{_F}.$

Now, by using the notation from \cite{ALM-Book},
we will introduce the notion of upper and lower functions,
that will be crucial to compute the rotation interval.

\begin{figure}
\centering
\begin{tikzpicture}[scale=0.25]
\draw (0,21.5) -- (0,0) -- (16,0) -- (16,21.5);
\draw (0,16) -- (16,16);
\node[below] at (0,0) {\tiny $0$}; \node[left] at (0,0) {\tiny $0$};
\node[below] at (16,0) {\tiny $1$}; \node[left] at (0,16) {\tiny $1$};
\begin{scope}[shift={(0,3)}]
   \draw[very thick] (-0.3,-1.5) -- (0,0) -- (4,9) -- (9,2) --
                     (12,11) -- (15,11) -- (16,16) -- (17, 18.25);
   \begin{scope}[shift = {(0.175,-0.175)}]
      \draw[red] (-0.3,-1.5) -- (0,0) -- (0.888,2) -- (9,2) -- (12,11) --
                (15,11) -- (16,16) -- (17, 18.25);
      \node[below] at (5,6.5) {$F$};
   \end{scope}
   \node[above left, blue] at (2.6,5.5) {$F_u$};
   \begin{scope}[shift={(-0.175,0.175)}]
      \draw[blue] (-0.3,-1.5) -- (0,0) -- (4,9) -- (11.3,9) -- (12,11) --
                  (15,11) -- (16,16) -- (17, 18.25);
   \end{scope}
   \node[below right, red] at (10.5,6.5) {$F_l$};
\end{scope}
\end{tikzpicture}
\caption{An example of a map $F \in \dol$ with its
lower map \textcolor{red}{$F_l$} in \textcolor{red}{red} and its
upper map \textcolor{blue}{$F_u$} in \textcolor{blue}{blue}.}\label{waterpicture}
\end{figure}
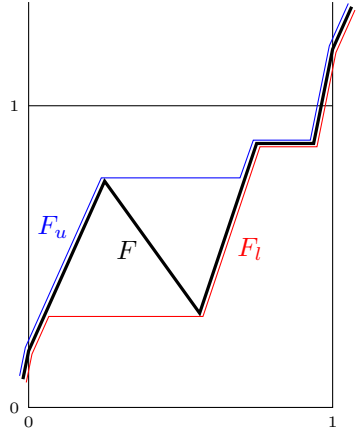

\begin{definition}\label{water}
Given $F\in\dol$ we define the \emph{$F$-upper map $F_u$} as
\[
F_u(x) := \sup\set{F(y)}{y\leq{}x}.
\]
Similarly we will define the \emph{$F$-lower map} as
\[
F_l(x) := \inf\set{F(y)}{y\geq{}x}.
\]
An example of such functions is shown in Figure~\ref{waterpicture}.
\end{definition}

It is easy to see that $F_l,F_u\in\dol$ are non decreasing, and
$F_l(x)\leq{}F(x)\leq{}F_u(x)$ for every $x\in\R.$

The rationale behind introducing the upper and lower functions comes
from the following result, stating that the rotation interval of a
function $F\in{}\dol$ is given by the rotation number of its
upper and lower functions.

\begin{theorem}\label{rotint}
Let $F\in\dol.$ Then,
\[
    \rot(F) = \bigl[\rho_{_{F_l}},\rho_{_{F_u}}\bigr].
\]
\end{theorem}
Note that this theorem makes indeed sense, since the upper and lower
functions are non-decreasing and by Proposition~\ref{nondecreasing}
they have a single well defined rotation number.

Let $\map{f}{\SI}$ and let $z\in{}\SI.$
The \emph{$f$-orbit of $z$} is defined to be the set
\[
    \Orb_f(z) := \{z, f(z),f^2(z),\dots{},f^{n}(z),\dots\}.
\]
We say that $z$ is an \emph{$n$-periodic point of $f$} if
$\Orb_f(z)$ has cardinality $n.$
Note that this is equivalent to $f^n(z) = z$ and
$f^k(z)\neq{}z$ for every $k<n.$
In this case the set $\Orb_f(z)$ will be called an
\emph{$n$-periodic orbit} (or, simply, a \emph{periodic orbit}).

If we have a periodic orbit of a circle map, a natural question
that might arise is how it behaves at a lifting level.
This motivates the introduction of the notion of a \emph{lifted cycle}.

Given a set $A\subset \R$ and $m \in \Z$ we will denote
$
   A+m := \set{x+m}{x\in{}A}.
$
Analogously, we set
\[
 A + \Z := \set{x+m}{x\in{}A,\ m\in\Z}.
\]

\begin{definition}
Let $\map{f}{\SI}{}$ be a continuous map and let $F$ be a lifting of $f.$
A set $P\subset{}\R$ is called a \emph{lifted cycle of $F$} if $e(P)$ is a
periodic orbit of $f.$ Observe that, then $P = P + \Z.$
The period of a lifted cycle is, by definition, the period of $e(P).$
Hence, when $e(P)$ is an $n$-periodic orbit of $f,$
$P$ is called an \emph{$n$-lifted cycle},
and every point $x\in{}P$ will be called an
\emph{$n$-periodic\!$\pmod{1}$ point of $F$}.
\end{definition}

The relation between lifted orbits and rotation numbers is
clarified by the next lemma.

\begin{lemma}\label{krotnum}
Let $F\in\dol.$
Then, $x$ is an $n$-periodic\!$\pmod{1}$ point of $F$
if and only if there exists $k \in \Z$ such that
$F^n(x) = x + k$ but $F^j(x) - x \notin \Z$ for $j=1,2,\dots,n-1.$
In this case,
\[
   \rho_{_F}(x) = \lim_{m\to\infty}\frac{F^m(x)-x}{m} = \frac{k}{n}.
\]
Moreover, let $P$ be a lifted $n$-cycle of $F.$
Every point $x \in P$ is an $n$-periodic\!$\pmod{1}$ point of $F,$
and the above number $k$ does not depend on $x.$
Hence, for every $x\in{}P$ we have
$\rho_{_F}(P) := \rho_{_F}(x) = \frac{k}{n}.$
\end{lemma}

Now we can revisit Proposition~\ref{nondecreasing}:
\setcounter{theorem}{2}
\begin{prop}
Let $F\in{}\dol$ be non-decreasing.
Then, for every $x\in{}\R$ the limit
\[
   \rho_{_F} := \lim_{n\to\infty}\frac{F^n(x)-x}{n}
\]
exists and is independent of $x.$
Moreover, $\rho_{_F}$ is rational if and only if $F$ has a lifted cycle.
\end{prop}
\setcounter{theorem}{7}

In the next two subsections we will survey on two known algorithms that have
been already used to compute rotation numbers of non-differentiable and
non-invertible liftings from $\dol.$
The first one (Algorithm~\ref{alg:ClassicAlgo}) stems automatically
from the definition of rotation number (Definition~\ref{defrotnum});
the other one (Algorithm~\ref{alg:SimosAlgo}) is due to
Sim\'o et al. \cite{simo}.

\subsection{Algorithm~\ref{alg:ClassicAlgo}:
the numerical algorithm to compute the rotation interval
that stems from the definition of rotation number}\label{subsec:alg:ClassicAlgo}
The first algorithm to compute $\rho_{_F}$ consists in using
Proposition~\ref{nondecreasing} and the following approximation,
for $n$ large enough in relation to the desired tolerance:\newline
\begin{minipage}{0.35\textwidth}\setlength\parindent{1.2em}
\begin{align*}
 \rho_{_F} &= \lim_{m\to\infty}\frac{F^m(x)-x}{m} \\
           &\approx \frac{F^n(x)-x}{n}\biggr\rvert_{x=0} \\
           &= \frac{F^n(0)}{n}.
\end{align*}

The implementation of the computation of this approximation to the rotation number
can be found in the side algorithm pseudocode.

Since the maps from $\dol$ are defined so that
\[ F(x) = F\evalat{[0,1]}\bigl(\decpart{x}\bigr) + \lfloor{x}\rfloor, \]
we need to evaluate the function
$\textproc{floor}(\cdot) = \lfloor{}\thinspace\cdot\thinspace\rfloor$ once per
\end{minipage}\hfill\begin{minipage}{0.6\textwidth}
\vspace*{-10pt}
\begin{algorithm}[H]
\caption{\strut\\ Direct Algorithm pseudocode}\label{alg:ClassicAlgo}
\begin{algorithmic}
\Procedure{Rotation\_Number}{$F$, \texttt{error}}
  \State $n \gets$ \textproc{ceil}$\Bigl(\tfrac{1}{\mbox{\texttt{error}}}\Bigr)$
  \State $x \gets 0$
  \State $k \gets 0$
  \For{$i \gets 1, n$}
       \State $x \gets F(x)$
       \State $s \gets \Call{floor}{x}$
       \State $k \gets k + s$\Comment{$k = \lfloor F^n(0) \rfloor$}
       \State $x \gets x - s$\Comment{$x = \decpart{F^n(0)} = F^n(0) - k$}
  \EndFor
  \State \textbf{return} $\tfrac{k+x}{n}$
\EndProcedure
\end{algorithmic}
\end{algorithm}
\end{minipage}

\vspace*{0.75ex}

\noindent
iterate. 
So, for clarity and efficiency, it seems advisable to split
$F^n(0)$ as $\decpart{F^n(0)} + \lfloor F^n(0) \rfloor.$
The next lemma clarifies the computation error as a function of the number of  iterates.
In particular it explicitly gives the necessary number of iterates, given a fixed tolerance.

For every non-decreasing lifting $F\in\dol,$ and every $n\in\N$ we set
(see Figure~\ref{desigualtats})
\[
\ell_F(n) := \min_{x\in\R} \left\lfloor F^n(x) - x\right\rfloor
           = \min_{x\in[0,1]} \left\lfloor F^n(x) - x\right\rfloor.
\]
The second equality holds because $F$ has degree 1,
and hence $\ell_F(n)$ is well defined.

\begin{lemma}\label{error}
For every non-decreasing lifting $F\in\dol$ and $n \in \N$ we have
\begin{enumerate}[(a)]
\item either $F^n(z) = z + \ell_F(n) + 1$ for some $z \in \R,$ or\newline
      $x + \ell_F(n) \leq F^n(x) < x + \ell_F(n) + 1$ for every $x\in \R;$
\item $\frac{\ell_F(n)}{n} \leq \rho_{_F} \leq \frac{\ell_F(n)+1}{n};$ and
\item $\abs{\rho_{_F} - \tfrac{F^n(x)-x}{n}} < \frac{1}{n}$
      for every $x \in \R.$
\end{enumerate}
\end{lemma}

\begin{proof}
We will prove the whole lemma by considering two alternative cases.
Assume first that $F^n(z) = z + \ell_F(n) + 1$
for some $z \in \R.$ Then (a) holds trivially, and
Proposition~\ref{nondecreasing} and Lemma~\ref{krotnum}
imply that $\rho_{_F} = \frac{\ell_F(n)+1}{n}.$
So, Statement~(b) also holds in this case.
Now observe that from the definition of $\ell_F(n)$
we have
\begin{equation}\label{definitionofell}
\ell_F(n) \leq \lfloor F^n(x) - x\rfloor \leq F^n(x) - x
\end{equation}
for every $x\in\R.$
Moreover, there exists $k = k(x) \in \Z$ such that
$x \in [z+k, z+k+1)$ and, since $F$ is non-decreasing, so is $F^n.$
Thus,
\begin{multline*}
F^n(x) - x \leq F^n(z+k+1) - x =
    F^n(z) + k + 1 - x = \\
    \ell_F(n) + 1 + (z + k + 1 - x) < \ell_F(n) + 2,
\end{multline*}
by Remark~\ref{iterats-de-grau-1}.
\enlargethispage{1ex}
Consequently,
\[
  \rho_{_F} -\frac{1}{n} =
  \frac{\ell_F(n)}{n} \leq
  \frac{F^n(x)-x}{n} < \rho_{_F} + \frac{1}{n};
\]
which proves (c) in this case.

\begin{minipage}{0.6\textwidth}
\setlength\parindent{1.2em}
Now we consider the case
\[ F^n(x) \neq x + \ell_F(n) + 1 \]
for every $x \in \R.$
In view of the definition of $\ell_F(n),$ we cannot have
\[ F^n(x) - x > \ell_F(n) + 1 \]
for every $x \in \R.$
Hence,
by the continuity of $F^n(x) - x$ and \eqref{definitionofell},
\begin{equation}\label{statementa}
  \ell_F(n) \leq F^n(x) - x < \ell_F(n) + 1
\end{equation}
for every $x \in \R.$ This proves (a).

Now we prove (b).
We consider the functions:
$x \longmapsto \ell_F(n) + x,$
$F^n,$ and
$x \longmapsto \ell_F(n) + 1 + x.$
They are all non-decreasing and, by Remark~\ref{iterats-de-grau-1},
they belong to $\dol.$
Hence, by Proposition~\ref{nondecreasing},
\cite[Lemma~3.7.19]{ALM-Book} and
\eqref{statementa},
\begin{multline*}
   \ell_F(n) = \rho_{_{x \mapsto \ell_F(n) + x}} \leq
       \rho_{_{F^n}} \leq \\
       \rho_{_{x \mapsto \ell_F(n) + 1 + x}} = \ell_F(n) + 1.
\end{multline*}
Consequently,
\[
   \frac{\ell_F(n)}{n} \leq \rho_{_{F}} = \frac{\rho_{_{F^n}}}{n} \leq \frac{\ell_F(n) + 1}{n},
\]
and (b) holds. Moreover, \eqref{statementa} is equivalent to
\[
   \frac{\ell_F(n)}{n} \leq \frac{F^n(x) - x}{n} \leq \frac{\ell_F(n) + 1}{n},
\]
which proves (c).
\end{minipage}\hfill\begin{minipage}{0.35\textwidth}
\vspace*{-3ex}
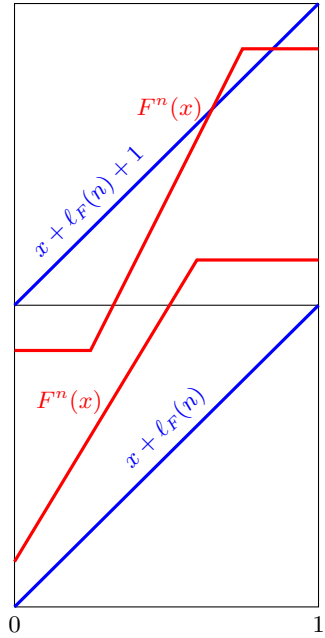
\begin{figure}[H]
\centering
\begin{tikzpicture}[scale=4]
\draw (0,0) -- (1,0) -- (1,2) -- (0,2) -- cycle; \draw (1, 1) -- (0, 1);
\draw[blue, very thick] (0,0) -- (1,1);
\draw[blue, very thick] (0,1) -- (1,2);

\draw[very thick, red] (0, 0.85) -- (0.25, 0.85) -- (0.75, 1.85) -- (1, 1.85);
\draw[very thick, red] (0, 0.15) -- (0.6, 1.15) -- (1, 1.15);

\node[below] at (0,0) {\small 0}; \node[below] at (1,0) {\small 1};

\node[left, red] at (0.65, 1.65) {\small $F^n(x)$};
\node[below left, red] at (0.32, 0.75) {\small $F^n(x)$};

\node[above right, rotate=45] at (0.1,1.1) {\small \textcolor{blue}{$x+\ell_F(n) + 1$}};
\node[above right, rotate=45] at (0.4,0.4) {\small \textcolor{blue}{$x+\ell_F(n)$}};
\end{tikzpicture}
\captionsetup{width=\linewidth}
\caption{Plot of $x+\ell_F(n)$ $x+\ell_F(n)+1,$
and $F^n(x)$ for two arbitrary non-decreasing maps $F\in\dol$ that fit in the two
cases of the lemma.}\label{desigualtats}
\end{figure}
\end{minipage}
\end{proof}

\subsection{Algorithm~\ref{alg:SimosAlgo}:
the Sim\'o et al. algorithm to compute the rotation interval}\label{subsec:alg:SimosAlgo}
First of all, it should be noted that even though the authors propose an algorithm to
compute the rotation interval for a general map $F\in\dol$, we will only use it for
non decreasing maps.
\enlargethispage{1ex}
\emph{A priori} this algorithm is radically different from Algorithm~\ref{alg:ClassicAlgo}
and it gives an estimate of $\rho_{_F}$ by providing and upper and a lower bound of the
rotation number (rotation interval in the original paper) of $F.$
Moreover, it is implicitly assumed that $\rho_{_F} \in [0,1]$
(in particular that $F(0) \in [0,1)$ --- this can be achieved by replacing the lifting
$F$ by the lifting $ G := F - \lfloor F(0) \rfloor,$ if necessary).
The algorithm goes as follows:
\begin{algorithm}[t]
\caption{Sim\'o et al. (\cite{simo}) Algorithm in pseudocode}\label{alg:SimosAlgo}
\begin{algorithmic}
\Procedure{Rotation\_Number}{$F,n$}
  \State \texttt{index[]} $\gets$
  \State $x \gets 0$
  \State $\rho_{\min} \gets 0$
  \State $\rho_{\max} \gets 1$
  \For{$i\gets 0, n$}
       \State $x \gets F(x)$
       \State $k_i \gets$ \Call{floor}{$x$}
       \State $\alpha_i \gets x-k_i$
       \State $\text{\ttfamily index[$i$]} \gets i$
  \EndFor
  \State \textbf{sort} $\alpha[\text{\ttfamily index[$i$]}]$ \textbf{by rearranging} \texttt{index[]}
  \For{$i\gets 0, n-1$}
       \State $\rho_{\mathrm{aux}} \gets \tfrac{k_{\text{\ttfamily index[$i+1$]}} - k_{\text{\ttfamily index[$i$]}}}{\text{\ttfamily index[$i+1$]} - \text{\ttfamily index[$i$]}}$
       \vspace*{4pt}
       \If{$\text{\ttfamily index[$i+1$]} > \text{\ttfamily index[$i$]}$}
           \State $\rho_{\min} \gets \max\{\rho_{\min}, \rho_{\mathrm{aux}}\}$
       \Else
           \State $\rho_{\max} \gets \min\{\rho_{\max}, \rho_{\mathrm{aux}}\}$
       \EndIf
  \EndFor
  \State \textbf{return} $\rho_{\min}, \rho_{\max}$
\EndProcedure
\end{algorithmic}
\end{algorithm}
\begin{enumerate}[({Alg. 2-}1)]
    \item Decide the number of iterates $n$ in function of a given tolerance.
    \item For $i = 0,1,2, \dots, n$ compute
          $k_i = \lfloor F^i(x_0) \rfloor$ and
          $\alpha_i = F^i(x_0) - k_i$ (i.e. $\alpha_i$ is the fractionary part of $F^i(x_0)$).
    \item Sort the values of $k_i$ and $\alpha_i$ so that
          $\alpha_{i_0} < \alpha_{i_1} < \dots < \alpha_{i_n}$
          (this can be achieved efficiently with the help of an index vector).
    \item Initialise $\rho_{\min} = 0$ and $\rho_{\max} = 1.$
    \item For $j = 0,1,2, \dots, n-1$ set $\rho_j = \frac{k_{i_{j+1}} - k_{i_{j}}}{i_{j+1}-i_{j}},$ and
    \begin{itemize}
       \item if $i_{j+1} > i_j$ set $\rho_{\min} = \max\{\rho_{\min}, \rho_j\};$ otherwise,
       \item if $i_{j+1} < i_j$ set $\rho_{\max} = \min\{\rho_{\max}, \rho_j\}.$
    \end{itemize}
    \item Return $\rho_{\max}$ and $\rho_{\min}$ as upper and lower bounds of the
          rotation number of $F,$ respectively.
\end{enumerate}

The real issue in this algorithm consists in dealing with the error.
If the rotation number $\rho_{_F}$ satisfies a Diophantine condition
$\abs{\rho_{_F} - \tfrac{p}{q}} \leq cq^{-\nu},$
with $c>0$ and $\nu \geq 2,$ then the error verifies
\[
   \varepsilon < \frac{1}{(c n^{\nu})^{\text{\tiny$\tfrac{1}{\nu-1}$}}}\,.
\]
Note that this error depends strongly on the chosen number $n$ of iterates,
and that $n$ \emph{must be chosen before knowing what the rotation number
could possibly be}.
Hence Algorithm~\ref{alg:SimosAlgo} it is not well suited to compute
\emph{unknown} rotation numbers of $\dol$ maps.
However, it is excellent in continuation methods where the current rotation
number gives a good estimate of the next one.

\begin{remark}
Note that the original aim of the algorithm to determine the existence of 
closed invariant curves on dynamical systems on the plane rather than the 
computation of rotation numbers of a given map of the circle.
The rationale of the algorithm is that if, 
after computing $\rho_{\min}$ and $\rho_{\max}$, 
we find that $\rho_{\min} > \rho_{\max}$ then the computed orbit 
cannot lay on a closed invariant curve.
This explains most of the limitations we have encountered, 
such as the lack of an \emph{a priori} estimate of the error, 
or the fact that the algorithm is suited only for rotation numbers $\rho\in[0,1]$.
\end{remark} 
\section{An algorithm to compute rotation numbers of non-decreasing maps with a constant section}\label{sec:methods}

The \emph{diameter of an interval $K$}
which, by definition is equal to the absolute value of the difference between their endpoints,
will be denoted as $\diam(K).$

A \emph{constant section of a lifting} $F$ of a circle map is a
closed non-degenerate
(i.e. different from a point or, equivalently, with non-empty interior,
or such that $\diam(K) > 0$)
subinterval $K$ of $\R$ such that $F\evalat{K}$ is constant.
In the special case when $F \in \dol,$ we have that
$F(x+1) = F(x) + 1 \ne F(x)$ for every $x\in\R.$
Hence, $\diam(K) < 1.$

The algorithm we propose is based on Lemma~\ref{error} but, specially,
on the following simple proposition which allows us to compute
\emph{exactly} the rotation number of a non-decreasing lifting from
$\dol$ that has a constant section, provided that
$F^n(K) \cap (K + \Z) \neq \emptyset.$
In this sense, Proposition~\ref{prop:fundamental} has a completely different
strategical aim than Algorithm~1 and Lemma~\ref{error}, which try to
\emph{(costly) estimate} the rotation number.

\begin{prop}\label{prop:fundamental}
Let $F\in{}\dol$ be non-decreasing and have a constant section $K.$
Assume that there exists $n \in \N$ such that
$F^n(K) \cap (K + \Z) \neq \emptyset,$
and that $n$ is minimal with this property.
Then, there exists $\xi \in\R$ such that
$F^n(K) = \{\xi\} \subset K+m$ with
$m = \lfloor \xi - \min K \rfloor \in \Z,$
$\xi$ is an $n$-periodic\!$\pmod{1}$ point of $F,$
and $\rho_{_F} = \frac{m}{n}.$
\end{prop}

\begin{proof}
Since $K$ is a constant section of $F,$
$F(K)$ contains a unique point, and hence there exists $\xi \in\R$
such that $F^n(K) = \{\xi\}.$
Then, the fact that $F^n(K) \cap (K + \Z) \neq \emptyset$ implies that
$\xi \in K + m$ with $m = \lfloor \xi - \min K \rfloor \in \Z.$

Set $\widetilde{\xi} := \xi - m \in K.$  Then,
$
\bigl\{F^n\bigl(\widetilde{\xi}\bigr)\bigr\} = F^n(K) = \bigl\{\widetilde{\xi} + m\bigr\}.
$
Moreover, the minimality of $n$ implies that
$F^j\bigl(\widetilde{\xi}\bigr) - \widetilde{\xi} \notin \Z$ for $j=1,2,\dots,n-1.$
So, Lemma~\ref{krotnum} tells us that $\widetilde{\xi}$ (and hence $\xi$) is an
$n$-periodic\!$\pmod{1}$ point of $F.$
Thus, $\rho_{_F} = \frac{m}{n}$ by Proposition~\ref{nondecreasing}.
\end{proof}

As already said, Proposition~\ref{prop:fundamental} is a tool to compute
\emph{exactly} the rotation numbers  of non-decreasing liftings
$F\in{}\dol$ which have a constant section  and have a lifted cycle
intersecting the constant section (and hence having rational rotation number).
In the next subsection we shall investigate how restrictive are these
conditions, when dealing with computation of rotation numbers.

\subsection{On the genericity of Proposition~\ref{prop:fundamental}}
First observe that the fact that Proposition~\ref{prop:fundamental}
only allows the computation of rotation numbers of non-decreasing liftings
$F\in{}\dol$ which have a constant section is not restrictive at all.
Indeed, if we want to compute rotation intervals of \emph{non-invertible}
continuous circle maps of degree one, Theorem~\ref{rotint} tells us that
this is exactly what we want.

Clearly, one of the real restrictions that cannot be overcome in the above method to
compute \emph{exact} rotation numbers is that it only works for maps having
a rational rotation number.

On the other hand, we also have the formal restriction that
Proposition~\ref{prop:fundamental} requires that the map $F$ has a lifted cycle
intersecting the constant section
(indeed this is a consequence of the condition $F^n(K) \cap (K + \Z) \neq \emptyset$).
A natural question is whether this restriction is just formal or it is a real one.
In the next example we will see that the restriction is not superfluous since
there exist maps which do not satisfy it.

Consequently, Proposition~\ref{prop:fundamental} is useless in
computing the rotation numbers of non-decreasing liftings in $\dol$
which have a constant section and either irrational rotation number
or rational rotation number but do not have any lifted cycle intersecting
the constant section.
The only reasonable solution to these problems is to use an iterative algorithm
to estimate the rotation number with a prescribed error, such as
Algorithm~\ref{alg:ClassicAlgo}, Algorithm~\ref{alg:SimosAlgo} or others.

\begin{example}
\emph{There exist non-decreasing liftings in $\dol$ which have a constant section
and rational rotation number but do not have any lifted cycle intersecting the constant section:}
\ Let $F \in \dol$ be the map such that
$F(x) = F\evalat{[0,1]}\bigl(\decpart{x}\bigr) + \lfloor{x}\rfloor$
for every $x \in \R,$ and let
\vspace*{-3ex}
\begin{figure}[H]
\raisebox{20ex}{\begin{minipage}{0.5\textwidth}
\[
   F\evalat{[0,1]}(x) := \begin{cases}
             x + 0.2            & \text{if $x \in [0,   0.1]$,}\\
             \frac{x}{2} + 0.25 & \text{if $x \in [0.1, 0.3]$,}\\
             7x - 1.7           & \text{if $x \in [0.3, 0.4]$,}\\
             \frac{x}{4} + 1    & \text{if $x \in [0.4, 0.8]$,}\\
             1.2                & \text{if $x \in [0.8, 1]$.}
     \end{cases}
\]
\end{minipage}}\hfill
\begin{tikzpicture}[scale=4]
\draw (0, 1.25) -- (0, 0) -- (1, 0) -- (1, 1.25);
\draw (0, 1) -- (1, 1);
\draw[dashed] (0, 0) -- (1, 1);
\draw[dashed] (0, 1) -- (0.25, 1.25);

\draw[blue, very thick] (-0.1, 0.2) -- (0, 0.2) -- (0.1, 0.3) -- (0.3, 0.4) --
                        (0.4, 1.1) -- (0.8, 1.2) -- (1, 1.2) -- (1.1, 1.3);

\foreach \x/\y in {0.1/0.3, 0.3/0.4, 0.4/1.1, 0.8/1.2}{ \draw[dotted] (\x,0) -- (\x,\y) -- (0,\y); }

\foreach \x in {0, 0.1, 0.3, 0.8, 1}{
     \draw (\x,0) -- (\x, -0.03);
     \node[below] at (\x, -0.03) {\small $\x$};
} ; \draw (0.4,0) -- (0.4, -0.03); \node[below right] at (0.34, -0.03) {\small $0.4$};

\foreach \y in {0, 0.3, 0.4, 1, 1.1, 1.2}{
     \draw (0,\y) -- (-0.03,\y);
     \node[left] at (-0.03,\y) {\small $\y$};
} \node[below left] at (0,0.2) {\small $0.2$};
\node[left, blue] at (0.35, 0.8) {\small $F\evalat{[0,1]}$};
\end{tikzpicture}\vspace*{-1ex}
\captionsetup{width=\linewidth}
\caption{Example of a non-decreasing lifting in $\dol$
with a constant section and rational rotation number
which does not verify the assumptions of
Proposition~\ref{prop:fundamental}.}\label{counter}
\end{figure}
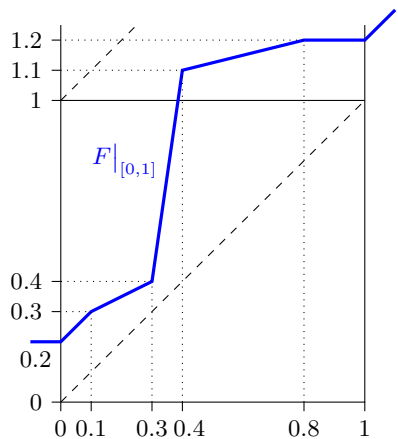

The map $F$ is a non-decreasing lifting from $\dol,$
having a constant section $K = [0.8, 1]$ and rotation number
$\tfrac{1}{3}$ given by the $3$-lifted cycle
$P = \{0.1, 0.3, 0.4\} + \Z$
(c.f. Lemma~\ref{krotnum} and Proposition~\ref{nondecreasing}).

Now let us see that $F$ does not have any lifted cycle
intersecting the constant section.
First, observe that
\[
  F^3(K) = F(F(F(K))) = F(F(\{1.2\})) = F(\{1.35\}) = \{1.75\} \not\subset K+\Z.
\]
Hence, there is no lifted cycle of period 3 intersecting $K.$
On the other hand, again by Lemma~\ref{krotnum}, we have that
if $x$ is an $n$-periodic\!$\pmod{1}$ point of $F$ then there exists
$k \in \Z$ such that $F^n(x) = x + k$
and
\[
   \frac{1}{3} = \rho_{_F} = \lim_{m\to\infty}\frac{F^m(x)-x}{m}
               = \rho_{_F}(x) = \frac{k}{n}.
\]
Moreover, since $F$ is non-decreasing, we know by
\cite[Corollary~3.7.6]{ALM-Book} that $n$ and $k$ must
be relatively prime. Thus, any lifted cycle of $F$ has period 3,
and from above this implies that there is no lifted cycle
intersecting $K.$
\end{example}
\subsection{Algorithm~\ref{alg:CSBAlgo}: A constant section based algorithm arising from Proposition~\ref{prop:fundamental}}
From the last paragraph of the previous subsection it becomes evident that
Proposition~\ref{prop:fundamental} does not give a complete algorithm to
compute rotation numbers of non-decreasing liftings in $\dol$
which have a constant section.
Such an algorithm must rather be a  mix-up of
Proposition~\ref{prop:fundamental}, and
Algorithm~\ref{alg:ClassicAlgo} to be used when we are not able to determine
whether we are in the assumptions of that proposition.
As in Algorithm~\ref{alg:ClassicAlgo},
for efficiency and because Proposition~\ref{prop:fundamental}
requires the computation of $m$ as an integer part,
we will split $F^n(0)$ as $\decpart{F^n(0)} + \lfloor F^n(0) \rfloor$
(here we are denoting the constant section by $K$ and assuming that $0\in K$
--- to be justified later).
Then, observe that the computations to be performed are exactly the same
in both cases
(meaning when we can use Proposition~\ref{prop:fundamental}, and when
alternatively we must end up by using Algorithm~\ref{alg:ClassicAlgo});
except for the conditionals that check whether there exists
$n \leq \mbox{\texttt{max\_iter}}$ such that
$F^n(K) \cap (K + \Z) \neq \emptyset$ is verified
(that is, whether the assumptions of Proposition~\ref{prop:fundamental} are verified)
before exhausting the \texttt{max\_iter} iterates determined a priori.

In what follows $\widetilde{F^n(0)}$ will denote the computed value of
$F^n(0)$ with rounding errors for $n=1,2,\dots,\mbox{\texttt{max\_iter}}.$

The algorithm goes as follows
(see Algorithm~\ref{alg:CSBAlgo} for a full implementation
in pseudocode, and see the explanatory comments below):
\begin{algorithm}[t]
\caption{Constant Section Based Algorithm\\
For a non-decreasing map $F \in \dol$ parametrised so that\\
$[-\mbox{\texttt{tol}}, \beta + \mbox{\texttt{tol}}]$
is a constant section of $F$}\label{alg:CSBAlgo}
\begin{algorithmic}
\State \textbf{define} \texttt{tol} $\gets$\Comment{\parbox[t]{20em}{%
                           Procedure parameter that bounds the rounding errors in the
                           computation of $F^n(0)$}}
\Statex
\Procedure{Rotation\_Number}{$F$, $\beta$, \texttt{error}}
  \State $\mbox{\texttt{max\_iter}} \gets \textproc{ceil}\Bigl(\tfrac{1}{\mbox{\texttt{error}}}\Bigr)$%
                          \Comment{\parbox[t]{17em}{Maximum number of iterates allowed
                              (to estimate the rotation number with the prescribed error when reached)}}
\vspace*{-4.5ex}
  \State $x \gets 0$
  \State $m \gets 0$
\vspace*{2pt}
  \For{$n \gets 1, \mbox{\texttt{max\_iter}}$}
     \State $x \gets F(x)$
     \State $s \gets \Call{floor}{x}$
     \State $m \gets m + s$\Comment{$m = \lfloor F^n(0) \rfloor$}
     \State $x \gets x - s$\Comment{$x = \decpart{F^n(0)} = F^n(0) - m$}
     \If{$x \leq \beta$}
         \State \textbf{return} $\tfrac{m}{n}$\Comment{\parbox[t]{24em}{%
                                      Exact rotation number:
                                      Proposition~\ref{prop:fundamental} holds assuming that the
                                      rounding error of $F^n(0)$
                                      is smaller than \texttt{tol}}}
\vspace*{-2ex}
     \EndIf
  \EndFor
   \State \textbf{return} $\tfrac{m+x}{\mbox{\texttt{max\_iter}}}$\Comment{\parbox[t]{20em}{%
                                      We do not know whether we are in the assumptions of
                                      Proposition~\ref{prop:fundamental}.
                                      So, we iteratively estimate the rotation number as in
                                      Algorithm~\ref{alg:ClassicAlgo}.
                                      The error bound is given by Lemma~\ref{error}}}
\vspace*{-7ex}
\EndProcedure
\vspace*{5ex}
\end{algorithmic}
\end{algorithm}
\begin{enumerate}[({Alg. 3-}1)]
    \item Decide the maximum number of iterates
    $\mbox{\texttt{max\_iter}} = \textproc{ceil}\Bigl(\tfrac{1}{\mbox{\texttt{error}}}\Bigr)$
    to perform in the worst case (i.e. when Proposition~\ref{prop:fundamental} does not work).
    \item Re-parametrize the lifting $F$ so that it has a maximal
         (with respect to the inclusion relation) constant section of the form
         $[-\mbox{\texttt{tol}}, \beta + \mbox{\texttt{tol}}],$
         where $\mbox{\texttt{tol}}$ is the pre-defined rounding error bound.
    \item Initialize $x = 0$ and $m = 0.$
    \item Compute iteratively
          $x = \decpart{\widetilde{F^n(0)}}$ and
          $m = \left\lfloor \widetilde{F^n(0)} \right\rfloor$
          (so that $\widetilde{F^n(0)} = x + m$)
          for $n \leq \mbox{\texttt{max\_iter}}.$
    \item Check whether $x \leq \beta$.
          On the affirmative we are in the assumptions of Proposition~\ref{prop:fundamental},
          and thus, $\rho_{_F} = \tfrac{m}{n}.$
          Then, the algorithm returns this value as the ``exact'' rotation number.
    \item If we reach the maximum number of iterates (i.e. $n = \mbox{\texttt{max\_iter}}$)
          without being in the assumptions of Proposition~\ref{prop:fundamental}
          (i.e. with $x > \beta$ for every $x$)
          then, by Lemma~\ref{error}, we have
\[
\abs{\rho_{_F} - \frac{m+x}{\mbox{\texttt{max\_iter}}}} =
\abs{\rho_{_F} - \frac{ \widetilde{F^n(0)}}{\mbox{\texttt{max\_iter}}}} \approx
\abs{\rho_{_F} - \frac{ F^n(0)}{\mbox{\texttt{max\_iter}}}} <
\frac{1}{\mbox{\texttt{max\_iter}}},
\]
and the algorithm returns $\frac{m+x}{\mbox{\texttt{max\_iter}}}$
as an estimate of $\rho_{_F}$ with $\tfrac{1}{\mbox{\texttt{max\_iter}}}$
as the estimated error bound.
\end{enumerate}

\begin{remark}
The fact that we can only \emph{check whether the assumptions of
Proposition~\ref{prop:fundamental} are verified before exhausting the}
$\mbox{\texttt{max\_iter}} = \textproc{ceil}\Bigl(\tfrac{1}{\mbox{\texttt{error}}}\Bigr)$
\emph{iterates determined a priori}
does not allow to take into account that $F$ may have
a lifted cycle intersecting the constant section but of very large period,
i.e. with period larger than \texttt{max\_iter}.
In practice this problem is totally equivalent to the non-existence
(or rather invisibility) of a lifted cycle intersecting the constant section,
and it can be considered as a new (algorithmic) restriction to
Proposition~\ref{prop:fundamental}.
It is solved in (Alg. 3-6) in the same manner as the two other problems related with the
applicability of Proposition~\ref{prop:fundamental} that have already been
discussed: by estimating the rotation number as in
Algorithm~\ref{alg:ClassicAlgo}.
\end{remark}

In the last part of this subsection we are going to discuss the rationale of
\hbox{(Alg.~3-2)} (and, as a consequence of (Alg.~3-5)).
The necessity of this tuning of the algorithm comes again from a challenge concerning
the application of Proposition~\ref{prop:fundamental}, which turns to be
one of the most relevant restrictions in the use of that proposition.
We will begin by discussing how we can efficiently check the condition
$\xi = F^n(0) \in K+\Z$  (or equivalently $F^n(K) \cap (K + \Z) \neq \emptyset$)
by taking into account that \emph{the computation of $F(x)$ is done with
rounding errors, and thus we do not know the exact values of $F^n(0)$ for}
$n=1,2,\dots,\mbox{\texttt{max\_iter}}$.
The next example shows the problems arising in this situation.

\begin{example}\label{example:Tangency}
\emph{$\widetilde{F^n(0)} \in K+\Z$ but $F^n(K) \cap (K + \Z) = \emptyset,$
and this leads to a completely wrong estimate of $\rho_{_F}$.}

\noindent Let $F \in \dol$ be the map such that
$F(x) = F\evalat{[0,1]}\bigl(\decpart{x}\bigr) + \lfloor{x}\rfloor$
for every $x \in \R,$ and let\newline
\raisebox{20ex}{\begin{minipage}{0.5\textwidth}
\[
F\evalat{[0,1]}(x) := \begin{cases}
    \tfrac{4}{3} x + \mu & \text{if $x \in \bigl[0, \tfrac{3}{4}\bigr],$}\\
                 1 + \mu & \text{if $x \in \bigl[\tfrac{3}{4}, 1\bigr],$}
\end{cases}
\]
with $\mu = \tfrac{819}{3124} - 10^{-16}.$
\end{minipage}}\hfill
\begin{tikzpicture}[scale=4]
\draw (0, 1.3) -- (0, 0) -- (1, 0) -- (1, 1.3);
\draw (0, 1) -- (1, 1);
\draw[dashed] (0, 0) -- (1, 1);
\draw[dashed] (0, 1) -- (0.3, 1.3);

\draw[blue, very thick] (0, 0.262) -- (0.75, 1.262) -- (1, 1.262);

\foreach \x in {0, 0.75, 1}{
     \draw (\x,0) -- (\x, -0.03);
     \node[below] at (\x, -0.03) {\small $\x$};
} \draw[dotted] (0.75, 0) -- (0.75, 1.262);

\foreach \y in {0, 0.262, 1, 1.262}{
     \draw (0,\y) -- (-0.03,\y);
     \node[left] at (-0.03,\y) {\small $\y$};
}
\node[left, blue] at (0.4, 0.85) {\small $F\evalat{[0,1]}$};
\end{tikzpicture}

For this map $F$ we have $K = \bigl[-\tfrac{3}{4}, 0\bigr]$ and
(see Figure~\ref{fig:GraphOfF5}) the graph of $F^5$ lies
above the graph of $x \longmapsto x+1$ and below
the graph of $x \longmapsto x+2,$ but very close to it
at five $F$-preimages of $x=\tfrac{3}{4}.$
On the other hand,
\[ F^5(0) = 1.74999999999999887\cdots \notin K + \Z \]
but the distance between $F^5(0)$ and $K + \Z$ is
$\tfrac{7}{4} - F^5(0) \approx 1.138\cdot 10^{-15}.$
Should the computations be done with rounding errors of this last
magnitude, we may have $\widetilde{F^5(0)} \gtrapprox\tfrac{7}{4},$
and accept erroneously that $F^5(0) \in K+\Z.$
This would lead to the conclusion that $\rho_{_F} = \tfrac{2}{5}$
but, as it can be checked numerically, $\rho_{_F} \approx 0.3983$
which is far from $\tfrac{2}{5}.$

At a first glance this seems to be paradoxical
but, indeed, it can be viewed in the following way:
The graph of $F^5$ does not intersect the
diagonal (modulo 1) $x+2,$ but there is a map $G$
close (at rounding errors distance) to $F$ such that
the graph of $G^5$ intersects that diagonal, and this
gives a lifted periodic orbit of period $5$ and rotation
number $\tfrac{2}{5}$ for $G.$
On the other hand, nothing is granted about the modulus
of continuity of $\rho_{_F}$ as a function of $F$
(notice that that everything here is continuous
including the dependence of the rotation number of $F$
on the parameter $\mu$), and this example explicitly shows
that it may be indeed very big.
In short, close functions can have very different rotation numbers.
\end{example}
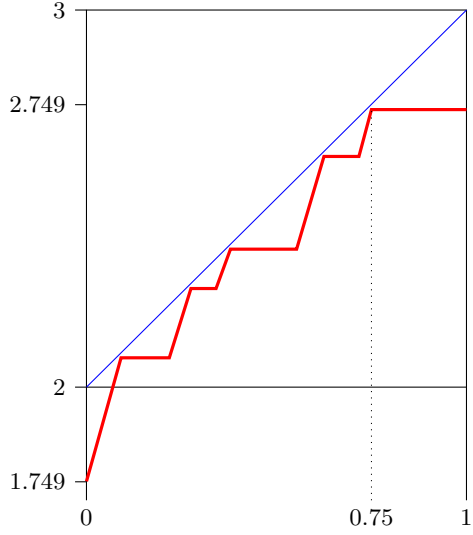
\begin{figure}[t]
\centering
\begin{tikzpicture}[scale=5]
\draw (0,1.7) -- (0,3) -- (1,3) -- (1,1.7); \draw (0,2) -- (1,2);
\draw[blue] (0,2) -- (1,3); 

\draw[very thick, red] (0, 1.749999) -- (0.091, 2.07778489) -- (0.218, 2.07778489) --
                       (0.275, 2.261473251) -- (0.341, 2.261473251) --
                       (0.379, 2.365877) -- (0.553, 2.365877) --
                       (0.625, 2.611715749) -- (0.717, 2.611715749) --
                       (0.75, 2.736) -- (1, 2.736);

\foreach \x in {0, 0.75, 1}{
     \node[below] at (\x, 1.7) {\small $\x$};
} \draw[dotted] (0.75, 1.7) -- (0.75, 2.736);

\foreach \y in {1.749, 2, 2.749, 3}{
     \draw (0,\y) -- (-0.03,\y);
     \node[left] at (-0.03,\y) {\small $\y$};
}
\end{tikzpicture}
\caption{The graph of $F^5.$ It lies below the graph of $x \longmapsto x+2$
but very close to it at five $F$-preimages of $x=\tfrac{3}{4}.$}\label{fig:GraphOfF5}
\end{figure}

The most reasonable solution to the problem pointed out in the
previous example consists in restricting the size of $K$ depending
of an a priori  estimate of the rounding errors in computing
$\widetilde{F^n(0)}$ for $n=1,2,\dots,\mbox{\texttt{max\_iter}}.$
Thus, we denote by \texttt{tol} an upper bound of these rounding errors, so that
\[
  \abs{F^n(0) - \widetilde{F^n(0)}} \leq \mbox{\texttt{tol}}
  \qquad
  \text{holds for } n=1,2,\dots,\mbox{\texttt{max\_iter}},
\]
and, given a maximal (with respect to the inclusion relation) constant section $K$
such that $0 \in K$ we write
$K := [\alpha - \mbox{ \texttt{tol}}, \beta + \mbox{ \texttt{tol}}].$
Then observe that the condition
$\widetilde{F^n(0)} \in [\alpha, \beta] + m$ for some $n\in\N$ and $m \in \Z$
implies $\xi = F^n(0) \in K + m,$
and $\rho_{_F} = \frac{m}{n}$ by Proposition~\ref{prop:fundamental}.

In practice, this ``rounding errors free'' version of the algorithm imposes
a new restriction to the applicability of Proposition~\ref{prop:fundamental}
(in the sense that it reduces even more the class of functions for which we
can get the ``exact rotation number'').
However, as before, the rotation numbers of the maps in the assumptions of
Proposition~\ref{prop:fundamental} for which we cannot compute
the ``exact rotation number'' can be estimated as in Algorithm~\ref{alg:ClassicAlgo}.

The computational efficiency of the algorithm strongly depends on how we check
the condition $\widetilde{F^n(0)} \in K+\Z.$
Taking into account the above considerations and improvements of the algorithm,
this amounts checking whether
$\alpha+\ell \leq \widetilde{F^n(0)} \leq \beta+\ell$ for some $\ell \in \Z,$
and we have to do so by using
$x = \decpart{\widetilde{F^n(0)}}$ and
$m = \left\lfloor \widetilde{F^n(0)} \right\rfloor$
instead of $\widetilde{F^n(0)} = x + m,$
which is the algorithmic available information.
Checking whether
$\alpha+\ell \leq \widetilde{F^n(0)} \leq \beta+\ell$ for some $\ell \in \Z$
is problematic since it requires at least two comparisons, and moreover
in general $\ell \neq m$
(and thus we need some more computational effort to find the right value of $\ell$).
A very easy solution to this problem is to change the parametrization of $F$
so that $\alpha = 0.$
In this situation we have
\[
m = m+\alpha \leq \widetilde{F^n(0)}, m+\beta < m + 1
\]
because $\diam(K) < 1,$ and $m = \left\lfloor \widetilde{F^n(0)} \right\rfloor.$
Consequently,
$\alpha+\ell \leq \widetilde{F^n(0)} \leq \beta+\ell$ for some $\ell \in \Z$
is equivalent to
\[
  \ell = m
  \quad\text{and}\quad
  x \le \beta.
\]
Thus, by ``tuning'' $F$ so that $\alpha = 0$ we get that $\ell = m$ and we manage
to determine whether $\widetilde{F^n(0)} \in [\alpha, \beta] + m$ just with one comparison
($x \le \beta$).

To see that and how we can change the parametrization of $F$
(that is the point 0) so that $\alpha = 0$ consider the map
$G(x) := F(x+\alpha) - \alpha.$
Clearly, $F$ and $G$ are conjugate by the
rotation of angle  $\alpha$: $x \longmapsto x + \alpha.$
Then, obviously, $G$ is a non-decreasing map in $\dol,$ has a constant section
$[-\mbox{ \texttt{tol}}, \beta - \alpha + \mbox{ \texttt{tol}}],$
and $\rho_{_F} = \rho_{_G}.$
So, every lifting can be replaced by one of its re-parametrizations with the same
rotation number and constant section
$[-\mbox{ \texttt{tol}}, \beta + \mbox{ \texttt{tol}}],$
where $\beta < 1 - 2\mbox{ \texttt{tol}}.$

\section{Testing the Algorithm}\label{sec:examples}
In this section we will test the performance of Algorithm~\ref{alg:CSBAlgo}
by comparing it against
Algorithms~\ref{alg:ClassicAlgo}~and~\ref{alg:SimosAlgo}
when dealing with different usual computations concerning rotation intervals.
First we will compare the efficiency of the three algorithms in computing and
plotting Devil's Staircases. Afterwards we will plot rotation intervals
and Arnold tongues for two  bi-parametric families that mimic the standard map
family. In the latter two cases, we will try to compare our algorithm with Algorithms~\ref{alg:ClassicAlgo}~and~\ref{alg:SimosAlgo} whenever possible.

\subsection{Computing Devil's staircases}\label{subsec:Compare}
In this subsection we will perform the comparison of algorithms by
computing  and plotting the Devil's staircase for the parametric family
$\bigl\{F_\mu\bigr\}_{\mu \in [0,1]} \subset \dol$
defined as
\begin{definition}\label{Fmu}
\[ F_\mu(x)=F_\mu\evalat{[0,1]}\bigl(\decpart{x}\bigr)+\lfloor{}x\rfloor,\]
where (see Figure~\ref{mupicture})
\begin{equation}\label{fmu}
F_\mu\evalat{[0,1]} (x) = \begin{cases}
\frac{4}{3}x+\mu &\text{if }x\leq\frac{3}{4}\\
\mu+1 &\text{if }x>\frac{3}{4}
\end{cases}.
\end{equation}
\end{definition}

Before doing this we shall remind the notion of a Devil's Staircase, 
and why typically exist for such families.
To this end we will first recall and survey on the notion of
\emph{persistence of a rotation interval}.

\begin{definition}
Given a subclass $\mathcal{A}$ of $\dol$, we say that
$F\in\mathcal{A}$ has an $\mathcal{A}$-persistent rotation interval
if there exists a neighbourhood $U$ of $F$ in $\mathcal{A}$ such that
\[ \rot(G) = \rot(F) \]
for every $G\in{}U.$
\end{definition}

We can now state the \emph{Persistence~Theorem} (c.f. \cite{persistence}):

\begin{theorem}[Persistence Theorem]\label{persistent}
Let $\mathcal{A}$ be a subclass of $\dol.$
Then the following statements hold:
\begin{enumerate}[(a)]
\item The set of all maps with $\mathcal{A}$-persistent rotation
interval is open and dense in $\mathcal{A}$ (in the topology of
$\mathcal{A}$).
\item If $F$ has an $\mathcal{A}$-persistent rotation interval, then
$\rho_{_{F_l}}$ and $\rho_{_{F_u}}$ are rational.
\end{enumerate}
\end{theorem}

\begin{remark}
If we apply Theorem~\ref{persistent} to our family
$\bigl\{F_\mu\bigr\}_{\mu \in [0,1]}$
which verifies that the rotation number of $F_\mu$
exists for every $\mu \in [0,1],$
we have that the set of parameters
$\mu \in [0,1]$ for which we have irrational rotation number has measure $0.$
Furthermore, for any $\kappa\in\Q$ such that there
exists $\mu$ with $\rho_{_{F_\mu}} = \kappa$, there exists an interval
$[\alpha,\beta]\ni{}\mu$ such that for all $\eta\in{}[\alpha,\beta],$
$\rho_{_{F_\eta}} = \kappa.$
\end{remark}

The so-called Devil's staircase is the result of plotting the
rotation number as a function of the parameter $\mu.$
By Theorem~\ref{persistent} we have that this plot will have
constant sections for any rational rotation number,
hence the "Staircase" in the name.

To test the algorithms, a $\mu$-parametric grid computation
of $\rho_{_{F_{\mu}}}$ with $\mu$ ranging from 0 to 1 with a step of $10^{-5}$
has been done.
For Algorithms~\ref{alg:ClassicAlgo} and \ref{alg:CSBAlgo} the
\texttt{error} has been set to $10^{-6}$.
For Algorithm~\ref{alg:CSBAlgo} the \texttt{tolerance}
has been set to $10^{-10}$.
For Algorithm~\ref{alg:SimosAlgo} we have arbitrarily set the number of iterates to $1000$.

\begin{figure}[h]
\centering
\includegraphics[width=0.5\linewidth, angle=-90]{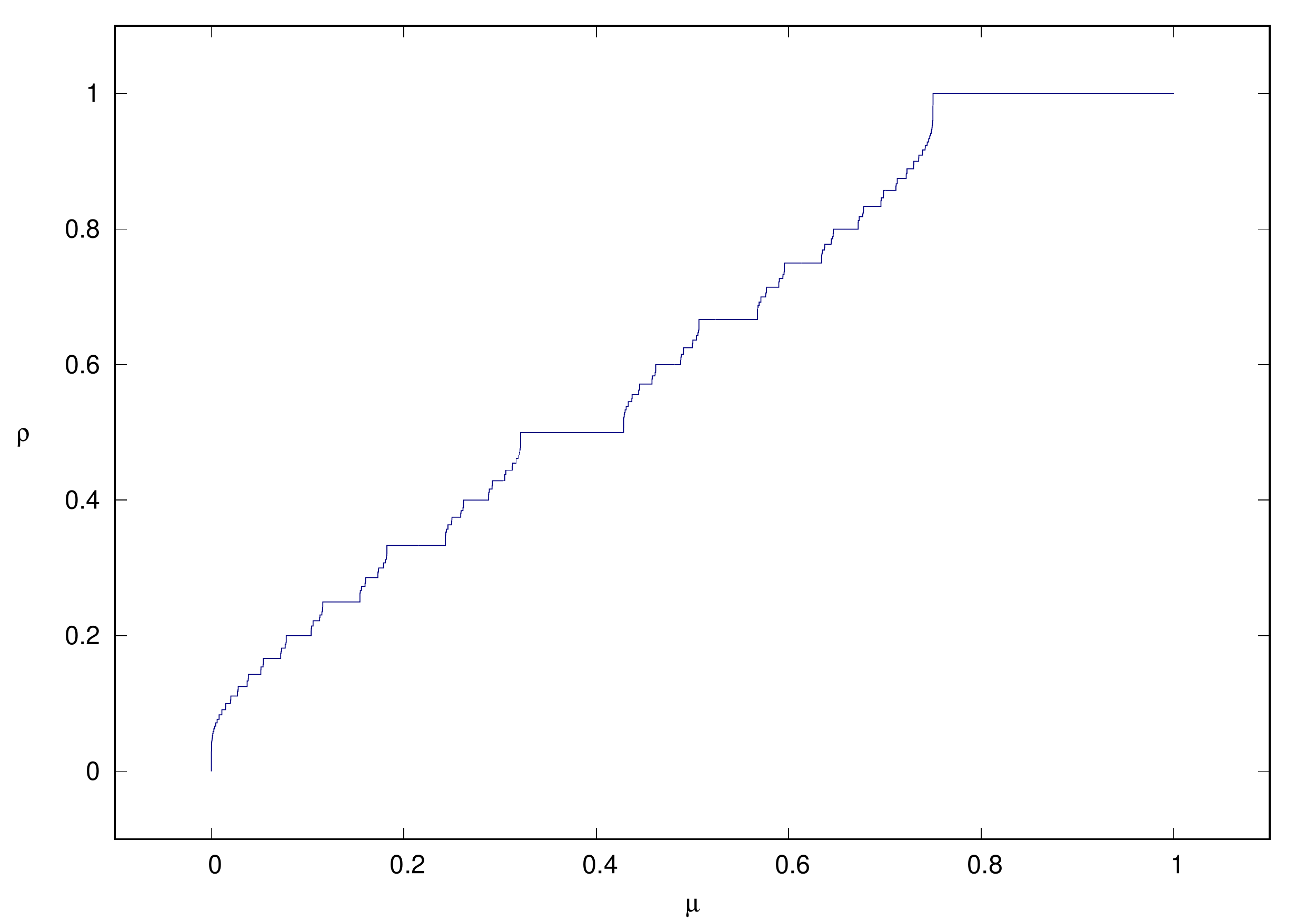}
\par\par\hfill
\includegraphics[width=0.275\linewidth, angle=-90]{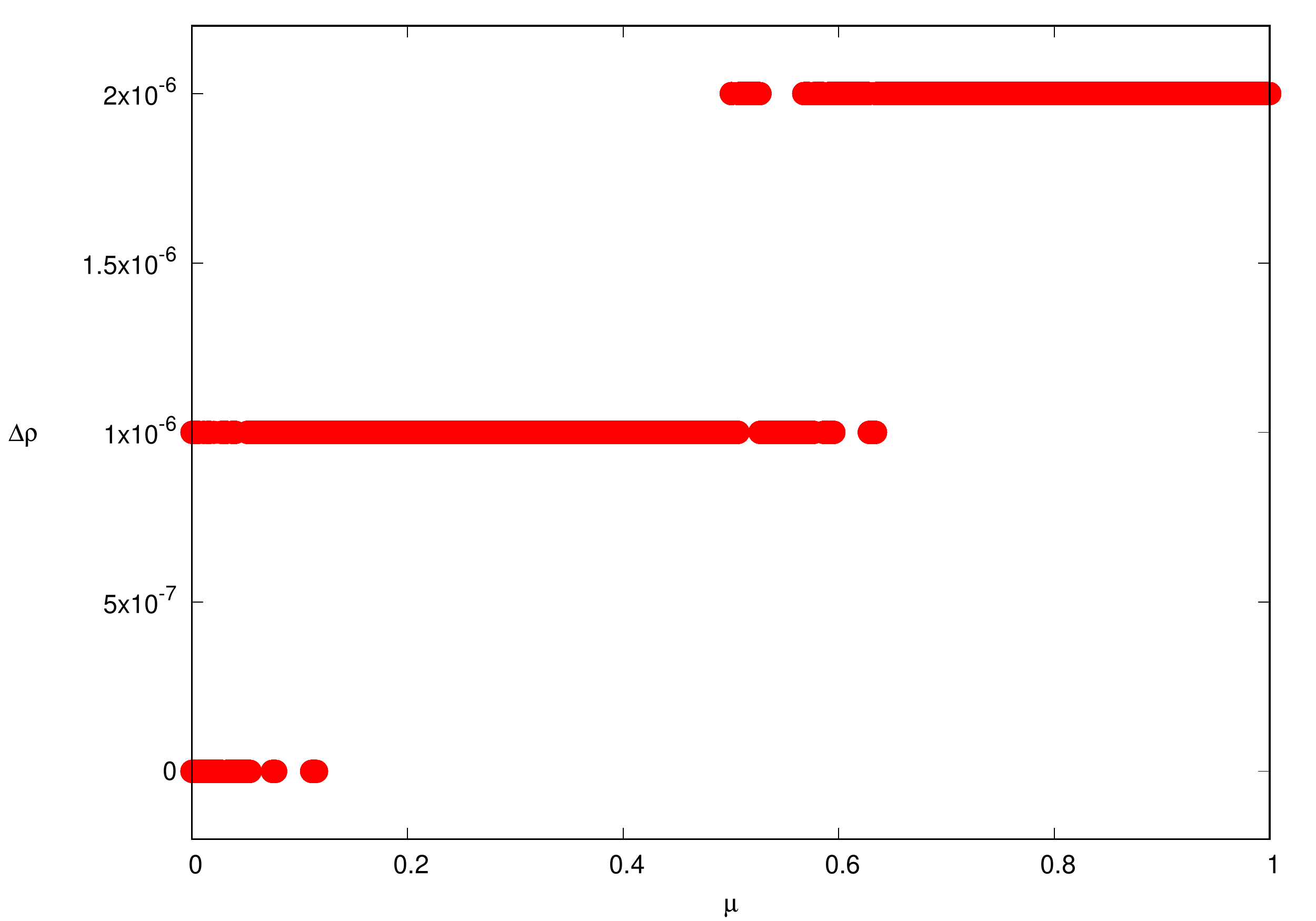}
\hfill
\includegraphics[width=0.275\linewidth, angle=-90]{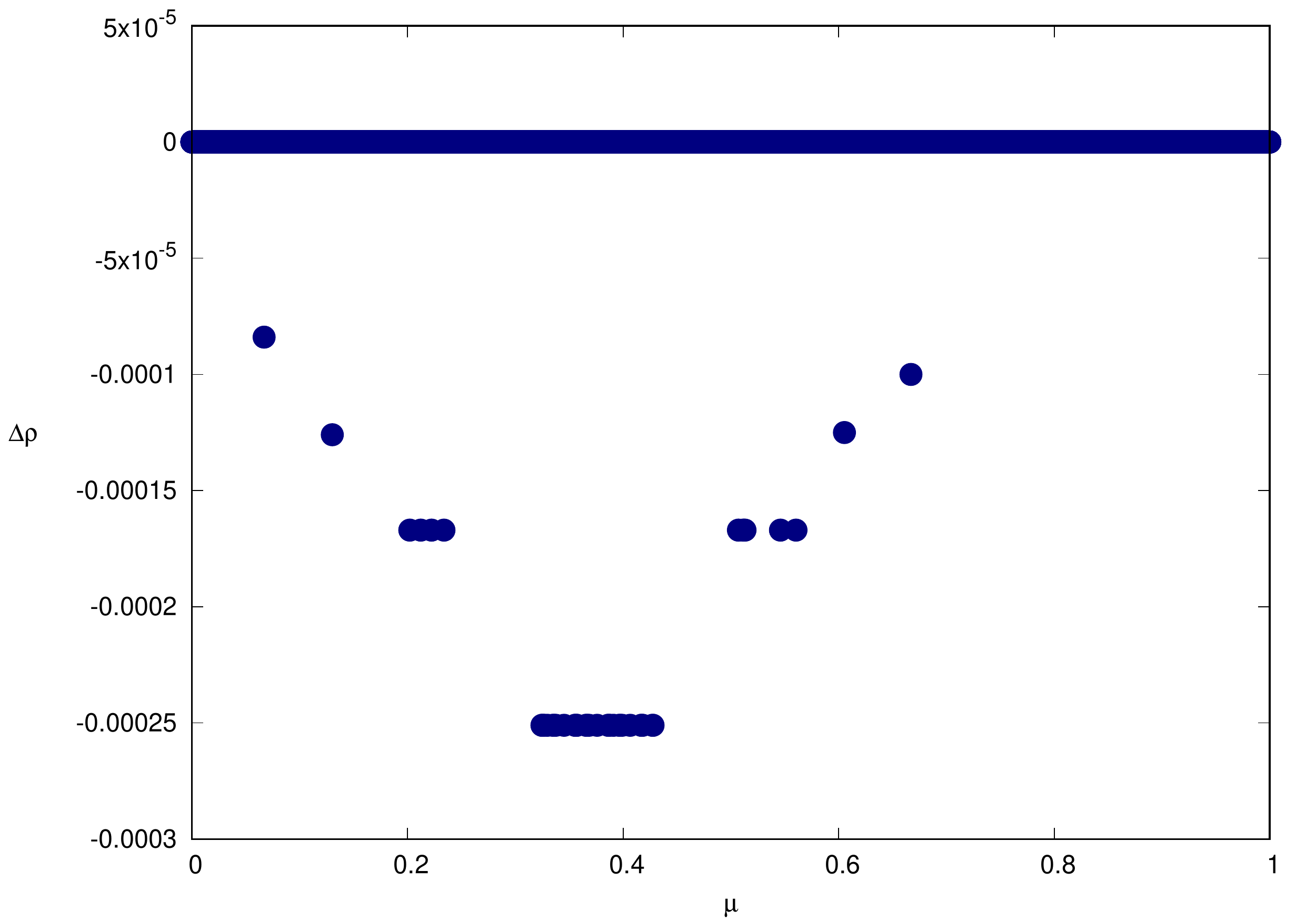}
\hfill\strut
\captionsetup{width=\linewidth}
\caption{The Devil's Staircase associated to the family \eqref{fmu}
computed with Algorithm~\ref{alg:CSBAlgo} (upper picture).
The lower pictures show the plots of the differences between
the value of $\rho_{_{F_{\mu}}}$ computed with Algorithm~\ref{alg:CSBAlgo} and
the value of $\rho_{_{F_{\mu}}}$ computed with Algorithm~\ref{alg:ClassicAlgo} (left picture), and
with the value of $\rho_{_{F_{\mu}}}$ computed with Algorithm~\ref{alg:SimosAlgo} (right picture).}\label{devils}
\end{figure}

In Figure~\ref{devils}
we show a plot of the Devil's Staircase computed with Algorithm~\ref{alg:CSBAlgo}, and
the plots of the differences between $\rho_{_{F_{\mu}}}$ computed with
Algorithms~\ref{alg:CSBAlgo}~and~\ref{alg:ClassicAlgo}, and the
differences between $\rho_{_{F_{\mu}}}$ computed with Algorithms~\ref{alg:CSBAlgo}~and~\ref{alg:SimosAlgo}.

Table~\ref{TimeTable}
\begin{table}
\begin{center}
\captionsetup{width=\linewidth}
\caption{Performance of the three algorithms studied for a variety of problems. The cells marked with N/A in \textcolor{blue}{blue} remark that Algorithm \ref{alg:SimosAlgo} does not work in general for $\rho\notin[0,1]$. The ones marked with N/A in \textcolor{red}{red} denote that the computation lasted more than a 100 processor hours and thus was terminated before it ended.}\label{TimeTable}
\begin{tabular}{llrcr}\toprule
\multirow{2}{*}{\textbf{Problem}} & \multirow{2}{*}{\textbf{Function Family}} & 
         \multicolumn{3}{c}{\textbf{Time taken by algorithm (s)}}\\ 
 & & \emph{Classic} & \emph{Sim\'o \emph{et al.}} &\emph{Proposed} \\ \midrule
Devil's Staircase & $F_\mu$ \small{(Def. \ref{Fmu})} & 2425.25 & 210.648 & 0.1413 \\ \midrule
\multirow{3}*{Rotation Interval}
           & Standard & 354.868 & \cellcolor{blue!25}N/A & 3.2874  \\
           & PWLSM \small{(Def.~\ref{Triangle})} & 110.892 & \cellcolor{blue!25}N/A & 0.4737\\
           & DSM \small{(Def. \ref{Disc})} & 63.588 & \cellcolor{blue!25}N/A & 0.2463\\ \midrule
\multirow{3}*{Arnol'd Tongues}
         & Standard & \cellcolor{red!25}N/A & \cellcolor{blue!25}N/A & 14948.41 \\
         & PWLSM & \cellcolor{red!25}N/A & \cellcolor{blue!25}N/A & 9729.17\\
         & DSM & \cellcolor{red!25}N/A & \cellcolor{blue!25}N/A & 4562.75\\
\bottomrule
\end{tabular}
\end{center}
\end{table}
shows the times\footnote{%
The simulations have been done with an Intel\textregistered{} Core\texttrademark{}
i7-3770 CPU @3.4GHz.}
taken by each of the three algorithms in computing the whole Devil's staircase using the three algorithms studied.

We remark that in the computation of the Devil's Staircase, Algorithm~\ref{alg:CSBAlgo} has been reduced to
Algorithm~\ref{alg:ClassicAlgo} only for  $\mu=0$ and for $\mu=1$,
as one would expect, since these cases follow the pattern of Example~\ref{example:Tangency}.

As a part of the testing of the algorithms we have also considered the inverse problem:
\emph{Given a value $x\in\R\backslash\Q$ and a tolerance $\varepsilon>0$
find the value $\mu = \mu(x)$ such that
$\rho_{F_\mu} \in [x-\varepsilon, x+\varepsilon].$}
This problem has turned to be extremely ill-conditioned:
by choosing $x$ to be an irrational such as the golden mean or $\pi/4$,
the continuity module of the function $\mu \mapsto \rho_{F_\mu}$
around $\mu(x)$ was estimated to be at least $10^{25}$,
making any attempt to solve the problem numerically a fool's errand.

\subsection{Rotation intervals for standard-like maps}\label{arnoldssect}
In this subsection we test our algorithm by efficiently computing the rotation
intervals and some \emph{Arnol'd tongues} for three bi-parametric families of maps:
the standard map family and two piecewise-linear extensions of it;
one continuous but not differentiable, and another one which is not even continuous.

We emphasize that the usual algorithms such as the ones from
\cite{simo2, Pavani, Seara-Villanueva, VanVeldhuizen}
cannot be used for these last two families families while the one we propose here
it works like a charm.

First we will recall the notion of \emph{Arnol'd tongue}.
\begin{definition}[Arnol'd Tongue \cite{Boyland}]
Let $\bigl\{F_{a,b}\bigr\}_{(a,b) \in \textsf{P}}$ be a two-parameter family of maps in $\dol$
for which the rotation interval $\rot\bigl(F_{a,b}\bigr)$ is well defined for every
possible point $(a,b) \in \textsf{P}$ in the parameter set.
Given a point $\varrho \in \R$ we define the $\varrho-$\emph{Arnold Tongue of}
$\bigl\{F_{a,b}\bigr\}_{(a,b) \in \textsf{P}}$ as
\[
\mathcal{T}_{\varrho} =
    \set{(a,b) \in \textsf{P}}{\varrho \in \rot\bigl(F_{a,b}\bigr)} \subset
    \textsf{P}.
\]
\end{definition}

Next we introduce each of the three families that we study and, for each of them we
show the results and we explain the performance of the algorithm.

\begin{figure}
\centering
\begin{tikzpicture}[scale=3]
\foreach \x in {-0.5, 0, 1, 1.5}{ \node[below] at (\x, -1) {\small \x}; \draw[dashed, gray] (\x,-1) -- (\x,2); }
\foreach \y in {-1, -0.5, 0, 0.5, 1, 1.5, 2}{ \node[left] at (-0.5, \y) {\small \y}; \draw[dashed, gray] (-0.5,\y) -- (1.5,\y); }
\draw[dashed, gray] (-0.5,-0.5) -- (1.5,1.5);
\draw[thick] (-0.5,-1) -- (1.5,-1) -- (1.5,2) -- (-0.5,2) -- cycle;

\draw[very thick] plot[smooth, samples=800, domain=-0.5:1.5] (\x, {\x - sin(360*\x)});

\draw[color=red, xshift=-0.5, yshift=0.5]
       plot[smooth, samples=200, domain=-0.5:-0.225] (\x, {\x - sin(360*\x)}) --
       (-0.225, 0.765) -- (0.5375, 0.765) --
       plot[smooth, samples=200, domain=0.537:0.775] (\x, {\x - sin(360*\x)}) --
       (0.775, 1.765) -- (1.52, 1.765);

\draw[color=blue, xshift=0.5, yshift=-0.5] (-0.52, -0.765) -- (0.225, -0.765) --
       plot[smooth, samples=200, domain=0.225:0.464] (\x, {\x - sin(360*\x)}) --
       (0.4635, 0.235) -- (1.225, 0.235) --
       plot[smooth, samples=200, domain=1.225:1.5] (\x, {\x - sin(360*\x)});
\end{tikzpicture}
\captionsetup{width=\linewidth}
\caption{The standard map with $a = 2\pi$ and $\Omega = 0$, with its
lower map in \textcolor{blue}{blue} and its
upper map in \textcolor{red}{red}.}\label{Standard}
\end{figure}
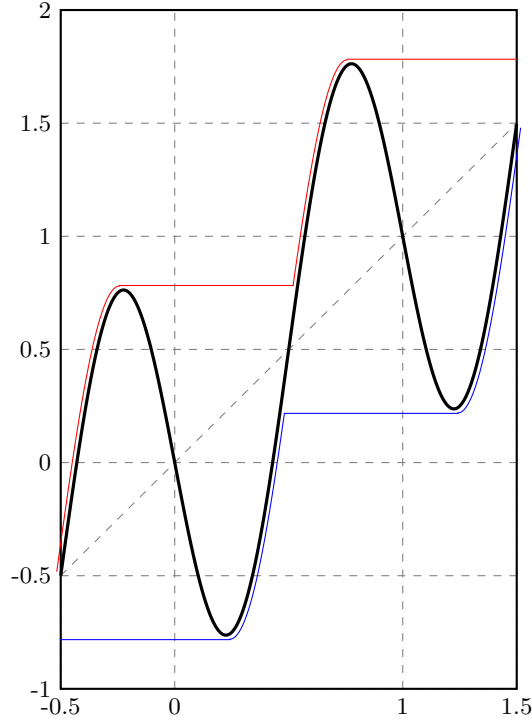

\bigskip\begin{definition}[Standard Map]
$\textsf{S}_{\Omega,a} \in \dol$ is defined as
(see Figure~\ref{Standard}):
\begin{equation}
\textsf{S}_{\Omega,a}(x) := x + \Omega - \frac{a}{2{}\pi} \sin(2\pi{}x).
\end{equation}
\end{definition}
To compute the rotation intervals of $\textsf{S}_{\Omega,a}$
we will use Theorem~\ref{rotint}, together with Algorithm~\ref{alg:CSBAlgo}.
To this end, first we will compute
$\bigl(\textsf{S}_{\Omega,a}\bigr)_l$ and $\bigl(\textsf{S}_{\Omega,a}\bigr)_u$
(that is, the lower and upper maps of $\textsf{S}_{\Omega,a}$), and then we will
use Algorithm~\ref{alg:CSBAlgo} to compute the rotation numbers
$\rho_{_{(\textsf{S}_{\Omega,a})_{_l}}}$ and
$\rho_{_{(\textsf{S}_{\Omega,a})_{_u}}}$ of these maps.

Note that $\textsf{S}_{\Omega,a}$ is non-invertible for $a > 1.$
Hence, in this case,
$\bigl(\textsf{S}_{\Omega,a}\bigr)_l$ and $\bigl(\textsf{S}_{\Omega,a}\bigr)_u$
do not coincide and have constant sections.
However, the characterization of these constants sections is not straightforward,
since their endpoints have to be computed numerically.
This is the reason why the computations of the rotation intervals and Arnol'd tongues
for the standard map have been the slowest ones.

\begin{figure}\centering
\includeepssubfigure{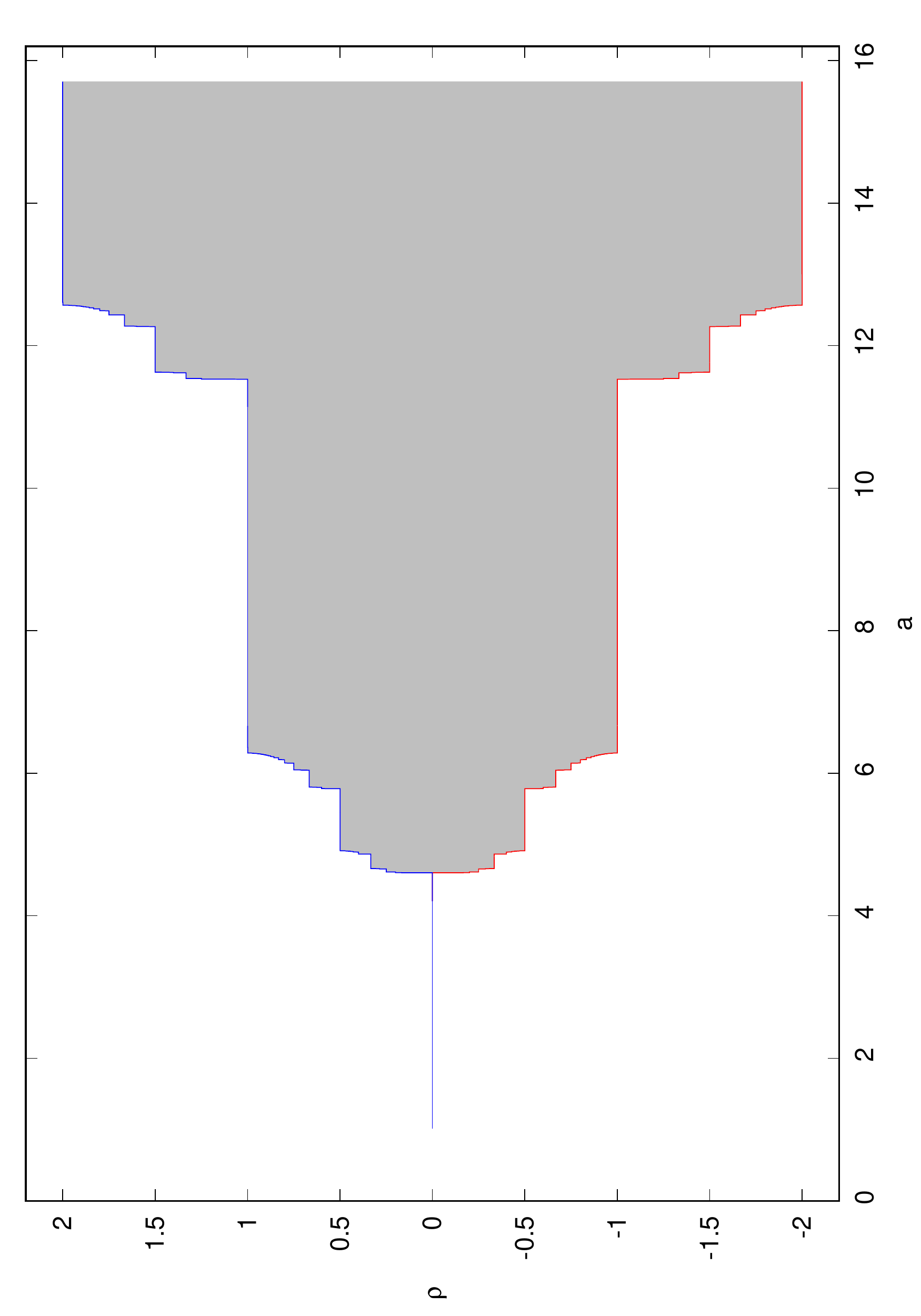}{Rotation interval graph for $\Omega = 0$.}
\hfill
\includeepssubfigure{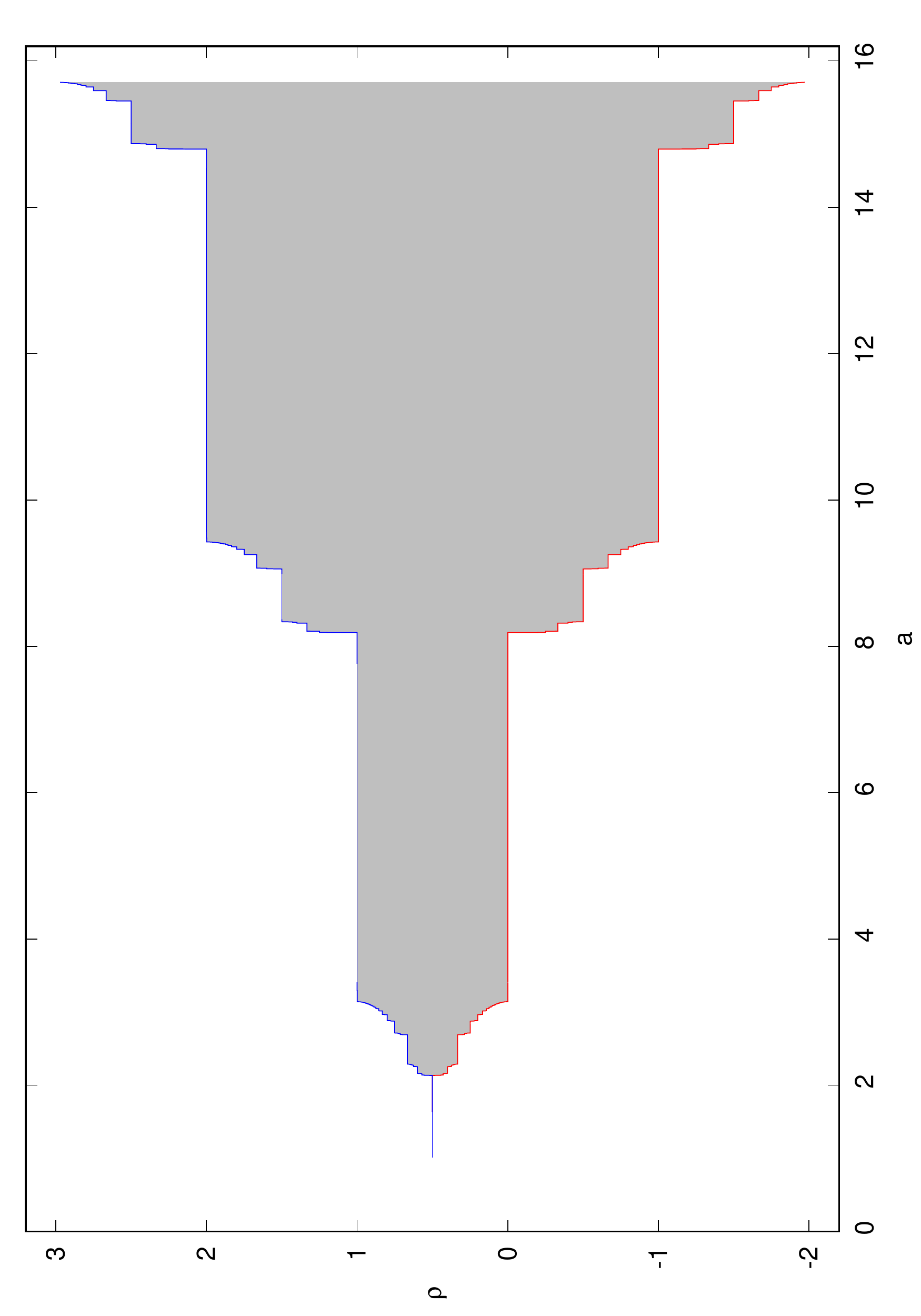}{Rotation interval graph for $\Omega = 0.5$.}
\hfill
\includeepssubfigure{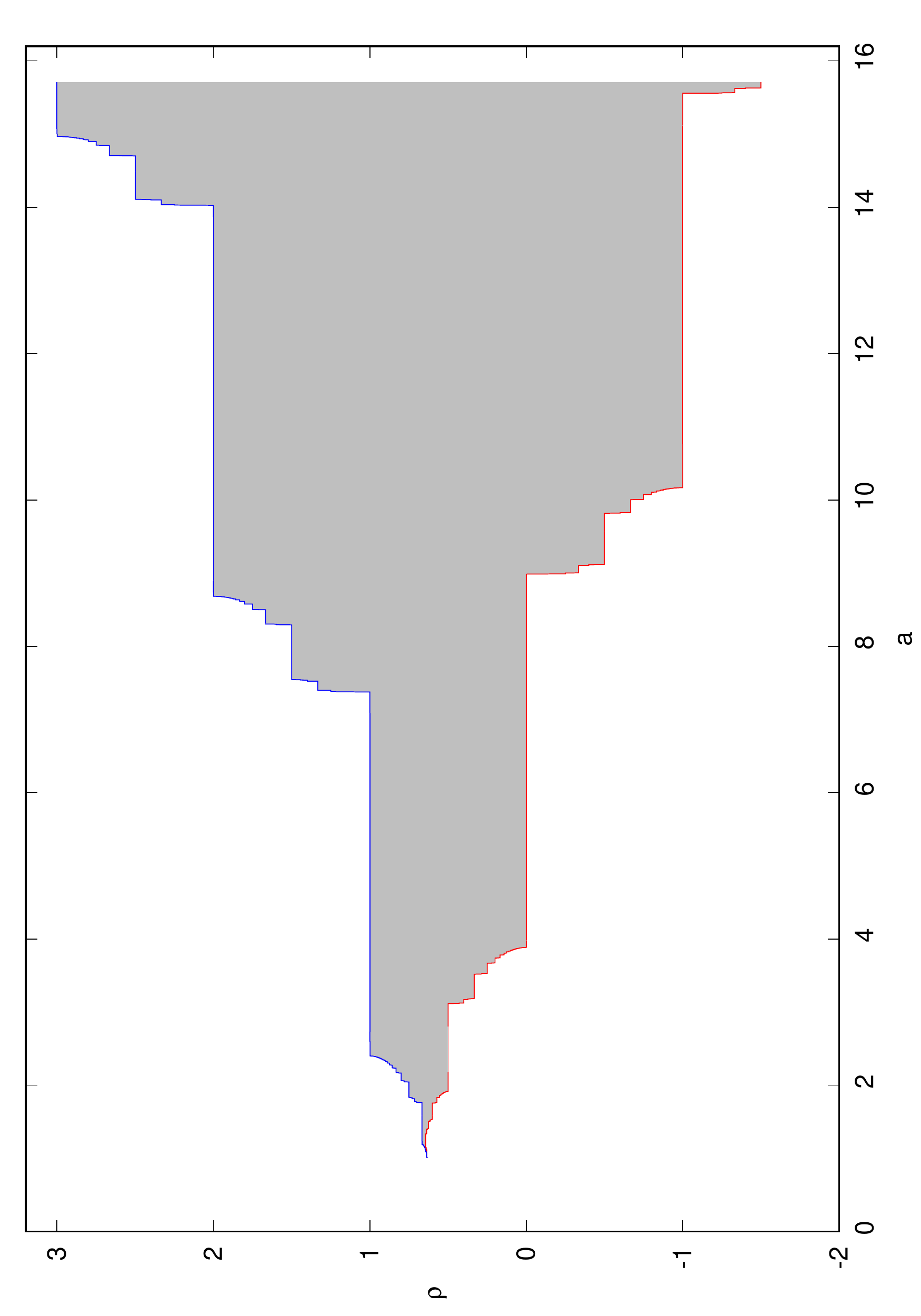}{Rotation interval graph with $\Omega$ equals to the Golden Mean.}
\par\par
\includepngsubfigure{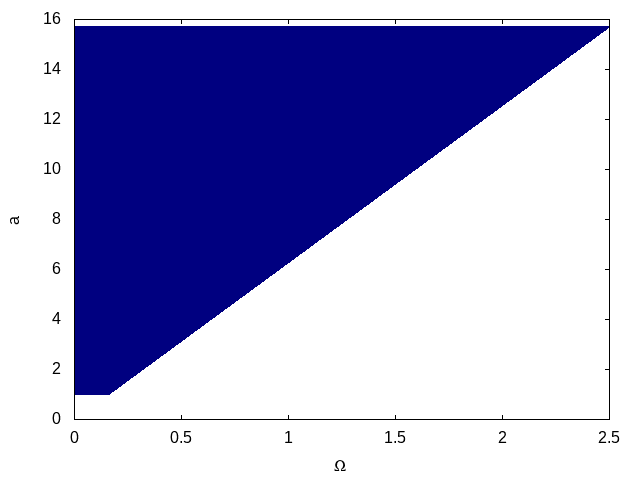}{$0-$Arnol'd ton\-gue.}
\hfill
\includepngsubfigure{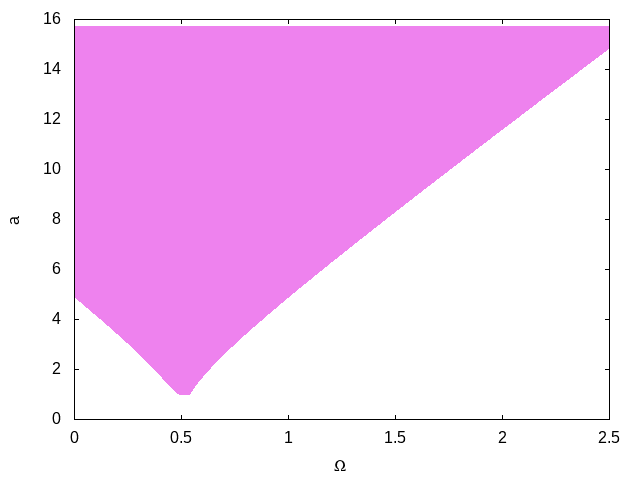}{$\mathcal{T}_{\varrho}$ Arnol'd tongue for $\varrho = 0.5$.}
\hfill
\includepngsubfigure{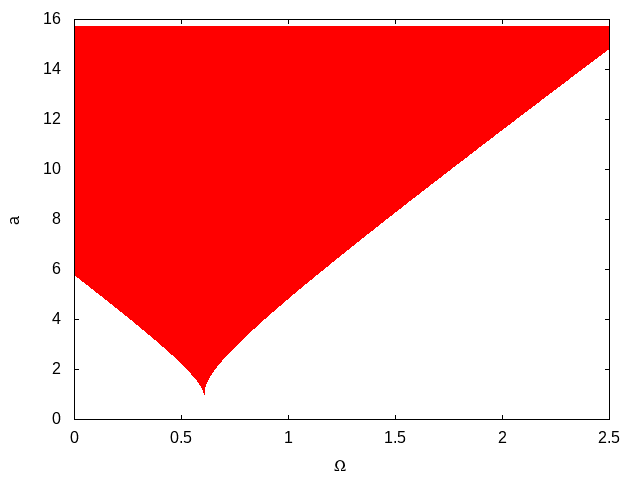}{$\mathcal{T}_{\varrho}$ for $\varrho$ equal to the Golden Mean.}
\captionsetup{width=\linewidth}
\caption{Graphs of the rotation interval and Arnol'd tongues for the Standard Map
$\textsf{S}_{\Omega,a}$.
The graphs of the rotation intervals are plotted as a function of the parameter $a$.}\label{GraphsStdMap}
\end{figure}

In Figure~\ref{GraphsStdMap} we show some graphs of the rotation interval and
Arnol'd tongues for the Standard Map. The graphs of the rotation intervals are
plotted for three different values of $\Omega$ as a function of the parameter $a$.

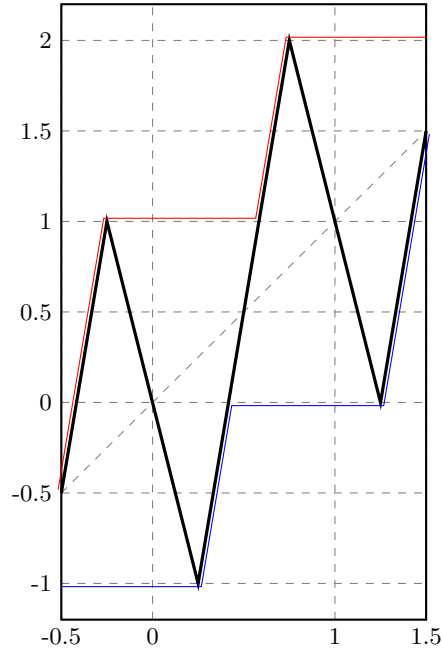
\begin{figure}[H]
\centering
\begin{tikzpicture}[scale=2.4]
\foreach \x in {-0.5, 0, 1, 1.5}{ \node[below] at (\x, -1.2) {\small \x}; \draw[dashed, gray] (\x,-1.2) -- (\x,2.2); }
\foreach \y in {-1, -0.5, 0, 0.5, 1, 1.5, 2}{ \node[left] at (-0.5, \y) {\small \y}; \draw[dashed, gray] (-0.5,\y) -- (1.5,\y); }
\draw[dashed, gray] (-0.5,-0.5) -- (1.5,1.5);
\draw[thick] (-0.5,-1.2) -- (1.5,-1.2) -- (1.5,2.2) -- (-0.5,2.2) -- cycle;

\draw[very thick] (-0.5,-0.5) -- (-0.25,1) -- (0.25,-1) -- (0.75, 2) -- (1.25,0) -- (1.5,1.5);
\draw[color=red, xshift=-0.5, yshift=0.5] (-0.5,-0.5) -- (-0.25,1) -- (0.583333333,1) -- (0.75, 2) -- (1.52,2);
\draw[color=blue, xshift=0.5, yshift=-0.5] (-0.52,-1) -- (0.25,-1) -- (0.4166666,0) -- (1.25,0) -- (1.5,1.5);
\end{tikzpicture}
\captionsetup{width=\linewidth}
\caption{The piecewise-linear standard map $\textsf{T}_{\Omega,a}$ with $a = \tfrac{5\pi}{2}$ and $\Omega = 0.$
The lower map of $\textsf{T}_{\Omega,a}$ is drawn in \textcolor{blue}{blue},
and the upper map in \textcolor{red}{red}.}\label{PLStandard}
\end{figure}

\begin{definition}[Piecewise-linear standard map]\label{Triangle}
We start by defining a convenience map\\[1pt]\noindent
\begin{minipage}{0.75\textwidth}
{\map{\tau}{[0,1]}[{[-1,1]}]} as follows:
\begin{equation}
\tau(x) = \begin{cases}
    4x      & \text{when $x \in \bigl[0, \tfrac{1}{4}\bigr]$,} \\
    2 - 4x  & \text{when $x \in \bigl[\tfrac{1}{4},\tfrac{3}{4}\bigr]$, and} \\
    4(x -1) & \text{when $x \in \bigl[\tfrac{3}{4},1\bigr]$.} \\
    \end{cases}
\end{equation}

Then, the \emph{piecewise-linear standard map} $\textsf{T}_{\Omega,a} \in \dol$ is defined by
(see Figure~\ref{PLStandard}):
\begin{equation}
\textsf{T}_{\Omega,a}(x) = x + \Omega - \frac{a}{2{}\pi}\tau\bigl(\decpart{x}\bigr),
\end{equation}
which corresponds to the standard map but using the $\tau$ wave function
instead of the $\sin(2\pi{}x)$ function.
\end{minipage}\hfill
\begin{minipage}{0.21\textwidth}\centering
\begin{tikzpicture}[scale=2]
\foreach \x in {0, 1,}{ \node[below] at (\x, -1) {\tiny \x}; }
\node[below] at (0.25, -1) {\tiny $\tfrac{1}{4}$}; \draw[densely dashed, gray] (0.25,-1) -- (0.25,1);
\node[below] at (0.75, -1) {\tiny $\tfrac{3}{4}$};
\foreach \y in {-1, 0, 1}{ \node[left] at (0, \y) {\tiny \y}; }
\draw (0, -1) -- (1,-1) -- (1,1) -- (0,1) -- cycle;
\draw (0, 0) -- (1,0);
\draw[dashed, gray] (0,-1) -- (1,0); \draw[dashed, gray] (0,0) -- (1,1);

\draw[blue, very thick] (0,0) -- (0.25,1) -- (0.75, -1) -- (1,0);
\end{tikzpicture}
\end{minipage}
\end{definition}
\noindent
The upper and lower maps for this family are very easy to compute.
Moreover, $\textsf{T}_{\Omega,a}$ is non-increasing for $a > \frac{\pi}{2}$ and hence,
in this case, the upper and lower maps do not coincide and have constant sections.

To compute the rotation intervals and Arnol'd Tongues of $\textsf{T}_{\Omega,a}$
we proceed as for the Standard Map by using Theorem~\ref{rotint} and
Algorithm~\ref{alg:CSBAlgo}.

\begin{figure}[H]
\includeepssubfigure{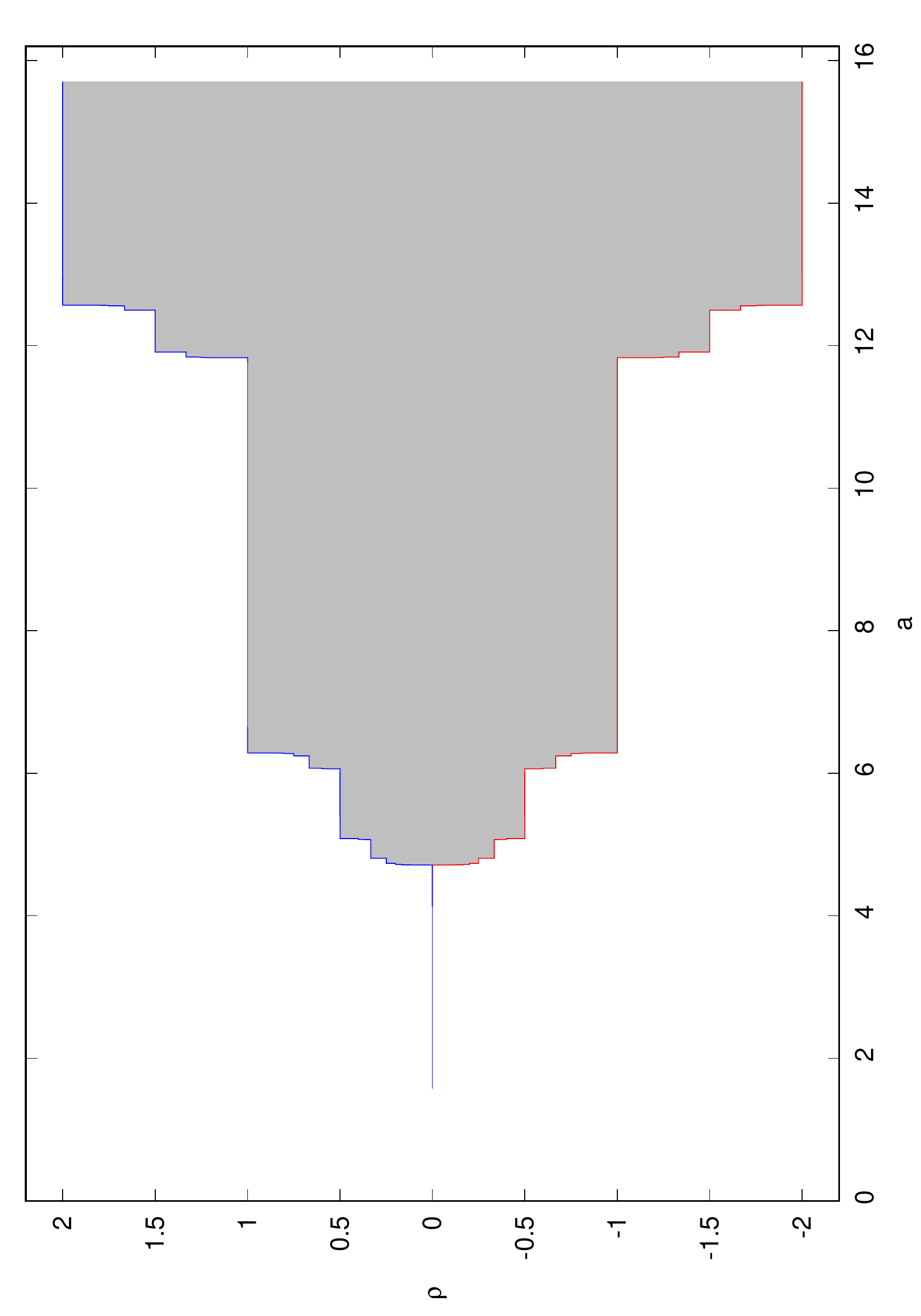}{Rotation interval graph for $\Omega = 0$.}
\hfill
\includeepssubfigure{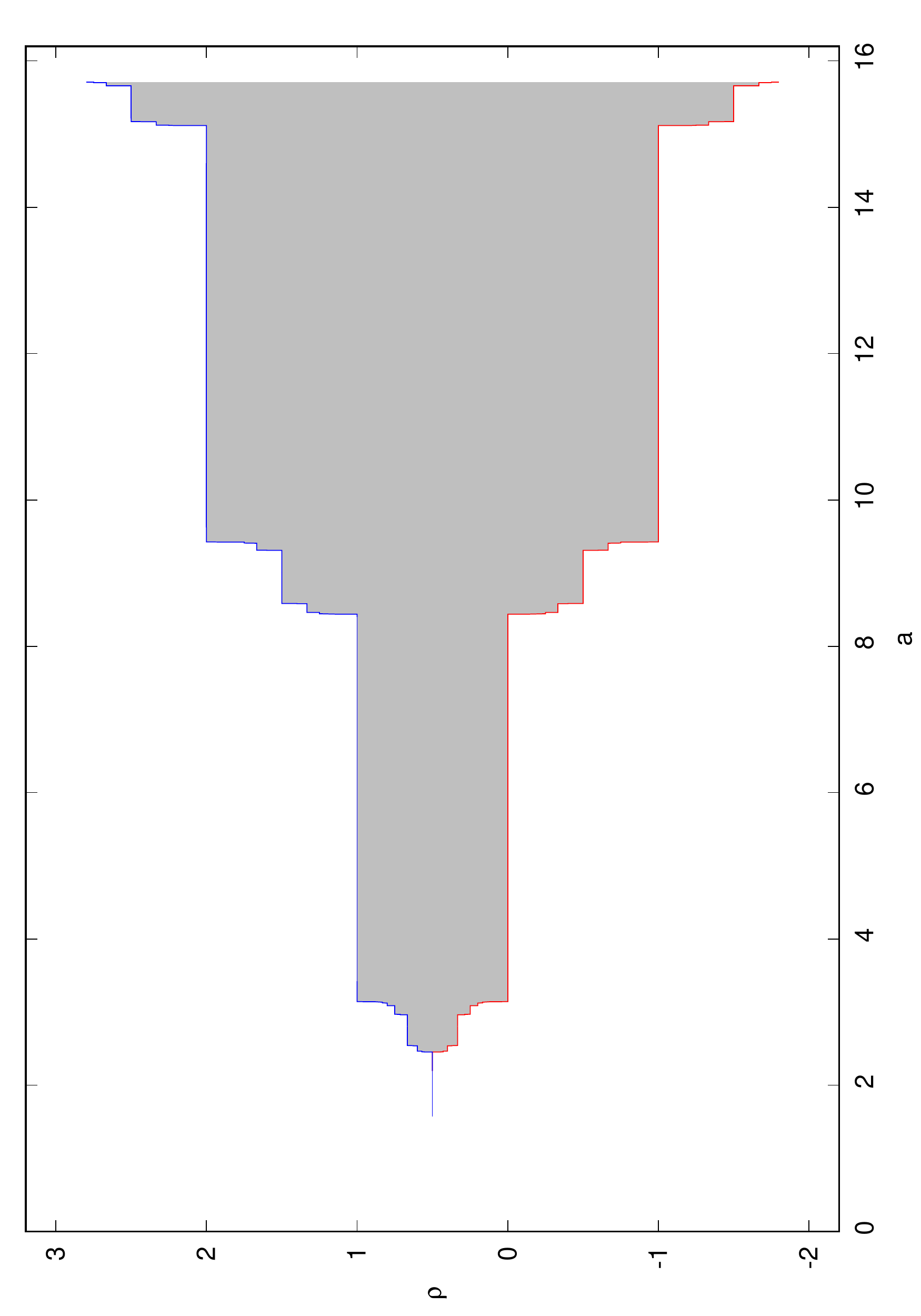}{Rotation interval graph for $\Omega = 0.5$.}
\hfill
\includeepssubfigure{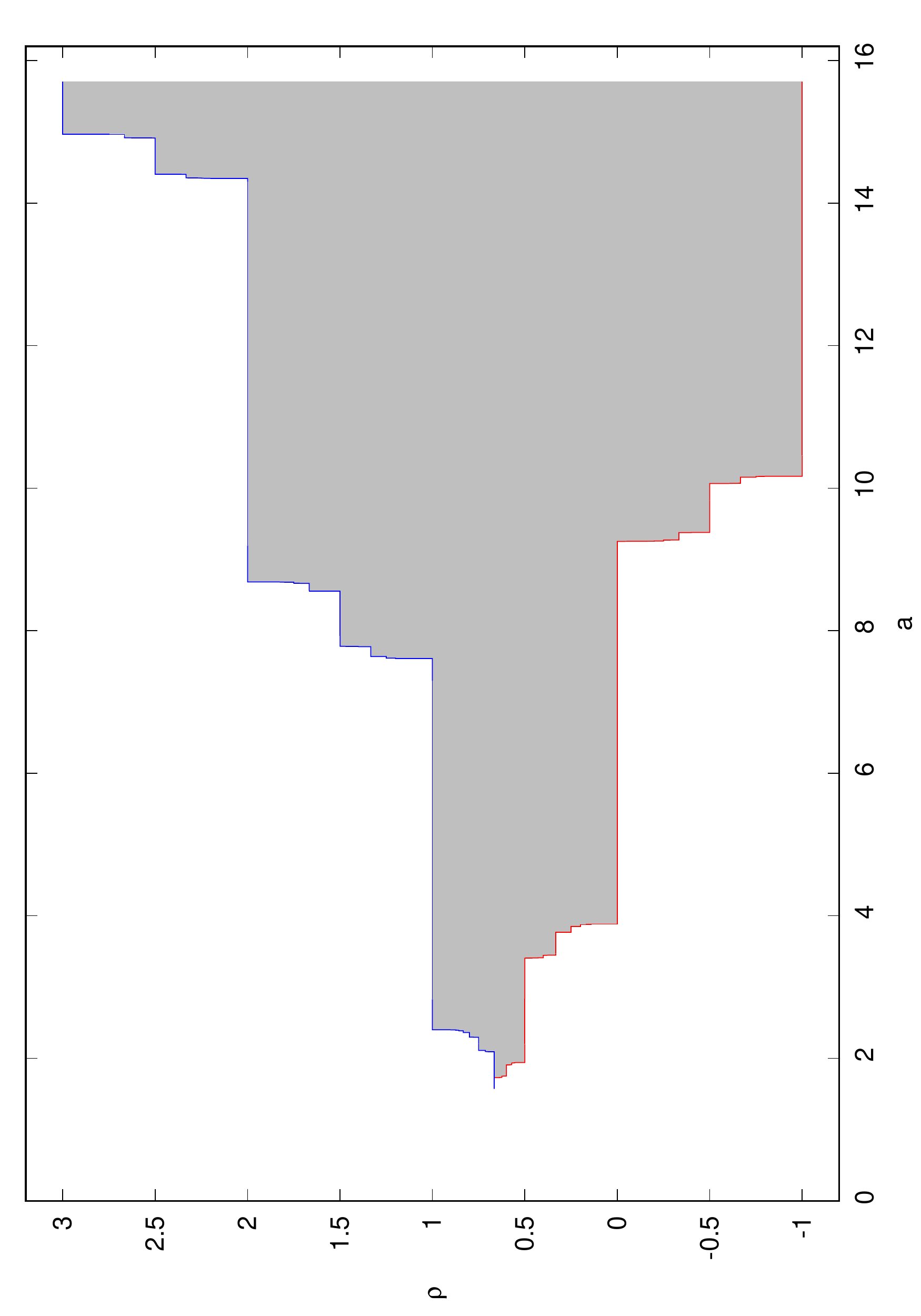}{Rotation interval graph with $\Omega$ equals to the Golden Mean.}
\par\par
\includepngsubfigure{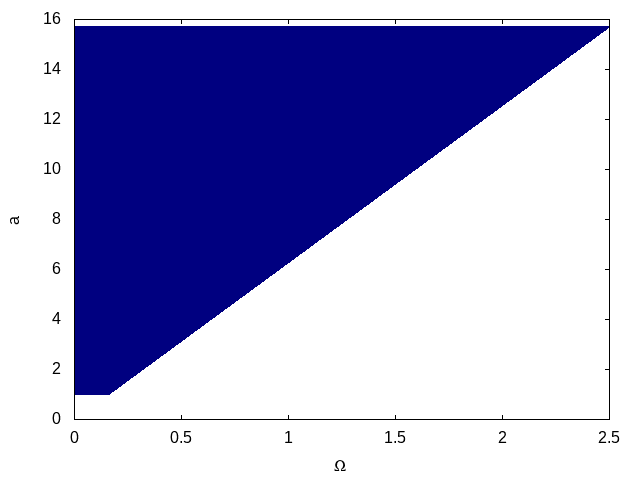}{$0-$Arnol'd ton\-gue.}
\hfill
\includepngsubfigure{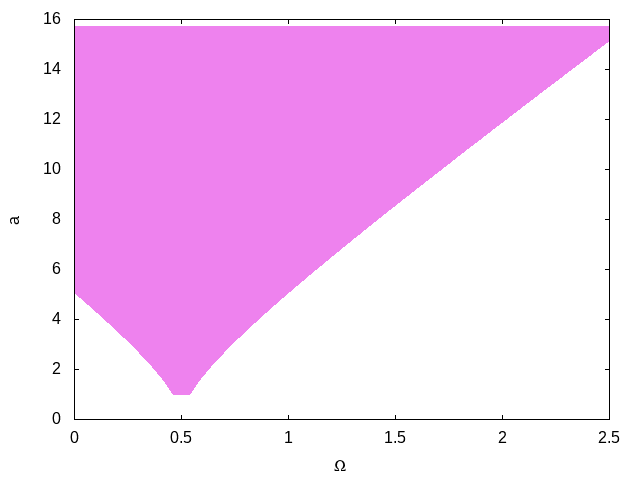}{$\mathcal{T}_{\varrho}$ Arnol'd tongue for $\varrho = 0.5$.}
\hfill
\includepngsubfigure{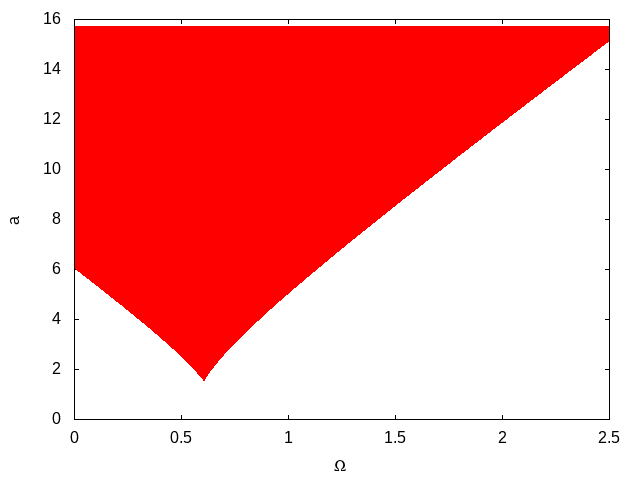}{$\mathcal{T}_{\varrho}$ for $\varrho$ equal to the Golden Mean.}
\captionsetup{width=\linewidth}
\caption{Graphs of the rotation interval and Arnol'd tongues for the piecewise-linear standard map
$\textsf{T}_{\Omega,a}$.
The graphs of the rotation intervals are plotted as a function of the parameter $a$.}\label{GraphsPLStdMap}
\end{figure}

In Figure~\ref{GraphsPLStdMap} we show some graphs of the rotation interval and
Arnol'd tongues for the piecewise-linear standard map.
The graphs of the rotation intervals are
plotted for three different values of $\Omega$ as a function of the parameter $a$.

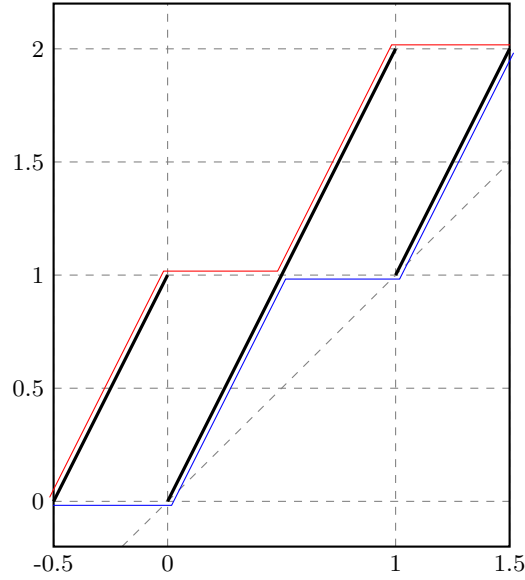
\begin{figure}[t]
\centering
\begin{tikzpicture}[scale=3]
\foreach \x in {-0.5, 0, 1, 1.5}{ \node[below] at (\x, -0.2) {\small \x}; \draw[dashed, gray] (\x,-0.2) -- (\x,2.2); }
\foreach \y in {0, 0.5, 1, 1.5, 2}{ \node[left] at (-0.5, \y) {\small \y}; \draw[dashed, gray] (-0.5,\y) -- (1.5,\y); }
\draw[dashed, gray] (-0.2,-0.2) -- (1.5,1.5);
\draw[thick] (-0.5,-0.2) -- (1.5,-0.2) -- (1.5,2.2) -- (-0.5,2.2) -- cycle;

\draw[very thick] (-0.5,0) -- (0,1);
\draw[very thick] (0,0) -- (1,2);
\draw[very thick] (1,1) -- (1.5,2);
\draw[color=red, xshift=-0.5, yshift=0.5] (-0.5,0) -- (0,1) -- (0.5,1) -- (1,2) -- (1.52,2);
\draw[color=blue, xshift=0.5, yshift=-0.5] (-0.52,0) -- (0,0) -- (0.5,1) -- (1,1) -- (1.5,2);
\end{tikzpicture}
\captionsetup{width=\linewidth}
\caption{The discontinuous standard map with $a = 2\pi$ and $\Omega = 0$ with its
lower map in \textcolor{blue}{blue} and its
upper map in \textcolor{red}{red}.}\label{DiscStandard}
\end{figure}

\bigskip\begin{definition}[The Discontinuous Standard Map]\label{Disc}
$\textsf{D}_{\Omega,a} \in \dol$ is defined as
(see Figure~\ref{DiscStandard}):
\begin{equation}
\textsf{D}_{\Omega,a}(x) := x + \Omega + \frac{a}{2{}\pi}\decpart{x}.
\end{equation}
\end{definition}

The map $\textsf{D}_{\Omega,a}$, being discontinuous, belongs to the so called 
class of \emph{old heavy maps} \cite{oldHeavy} 
(the \emph{old} part of the name stands for \emph{degree one lifting}
 --- that is, $\textsf{D}_{\Omega,a} \in \dol$).
A map $F \in \dol$ is called \emph{heavy} if for any $x \in \R$,
\[
  \lim_{y\searrow{}x^+} F(y) \leq F(x) \leq \lim_{y\nearrow{}x^-} F(y)
\]
(in other words, the map ``falls down'' at all discontinuities).

Observe that for the class of \emph{old heavy maps} the upper and lower maps
in the sense of Definition~\ref{water} are well defined and continuous.
Moreover, the whole family of \emph{water functions} (c.f. \cite{ALM-Book})
is well defined and continuous.
So, the rotation interval of the \emph{old heavy maps} is well defined
\cite[Theorem~A]{oldHeavy} and, moreover,
Theorem~\ref{rotint} together with Algorithm~\ref{alg:CSBAlgo} work for this class.
Hence, to compute the rotation intervals and Arnol'd Tongues of $\textsf{D}_{\Omega,a}$
we proceed again as for the Standard Map.

As for the piecewise-linear standard maps the upper and lower maps are very easy
to compute, and have constant sections for $a \neq 0$

\begin{figure}[h]
\includeepssubfigure{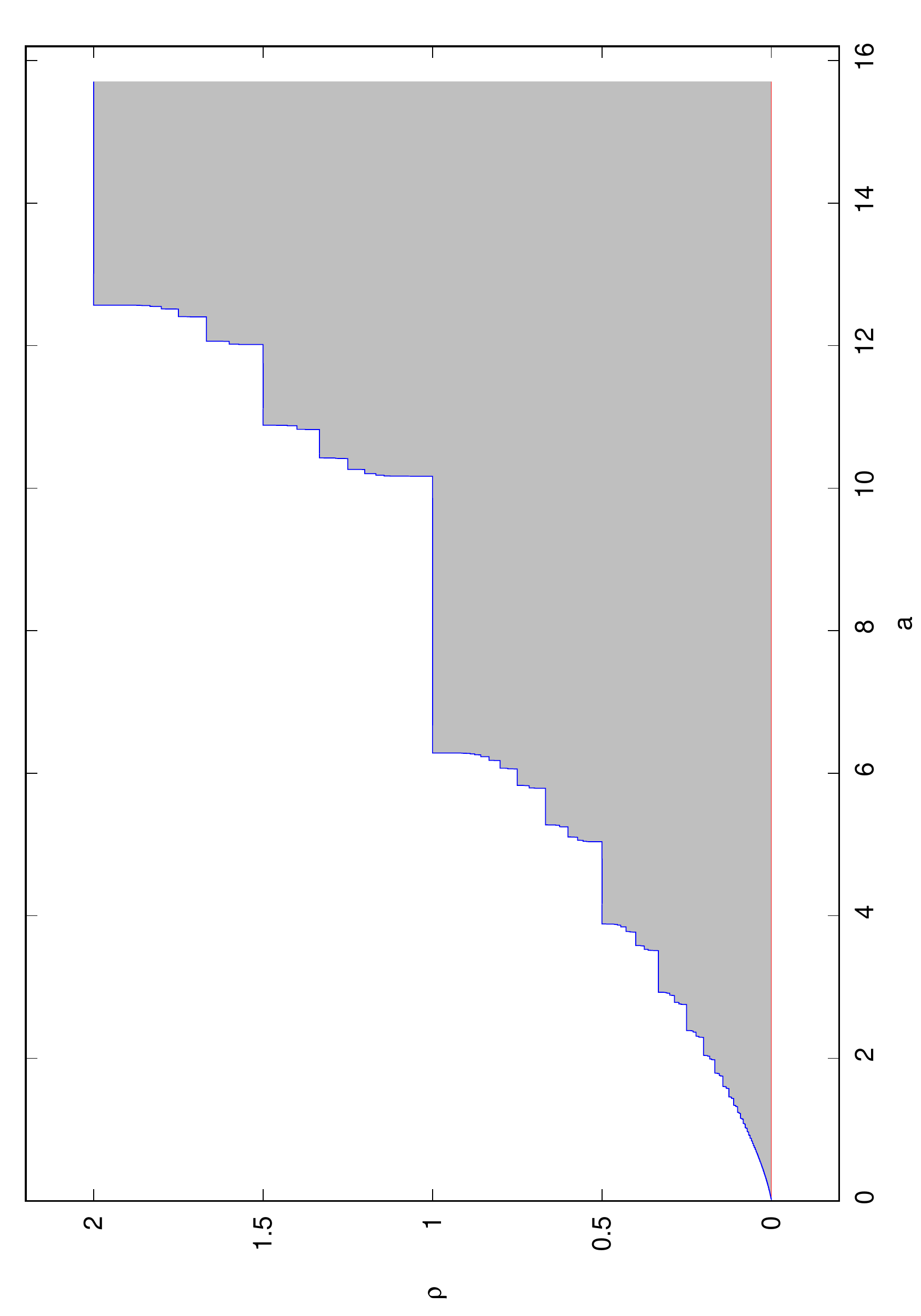}{Rotation interval graph for $\Omega = 0$.}
\hfill
\includeepssubfigure{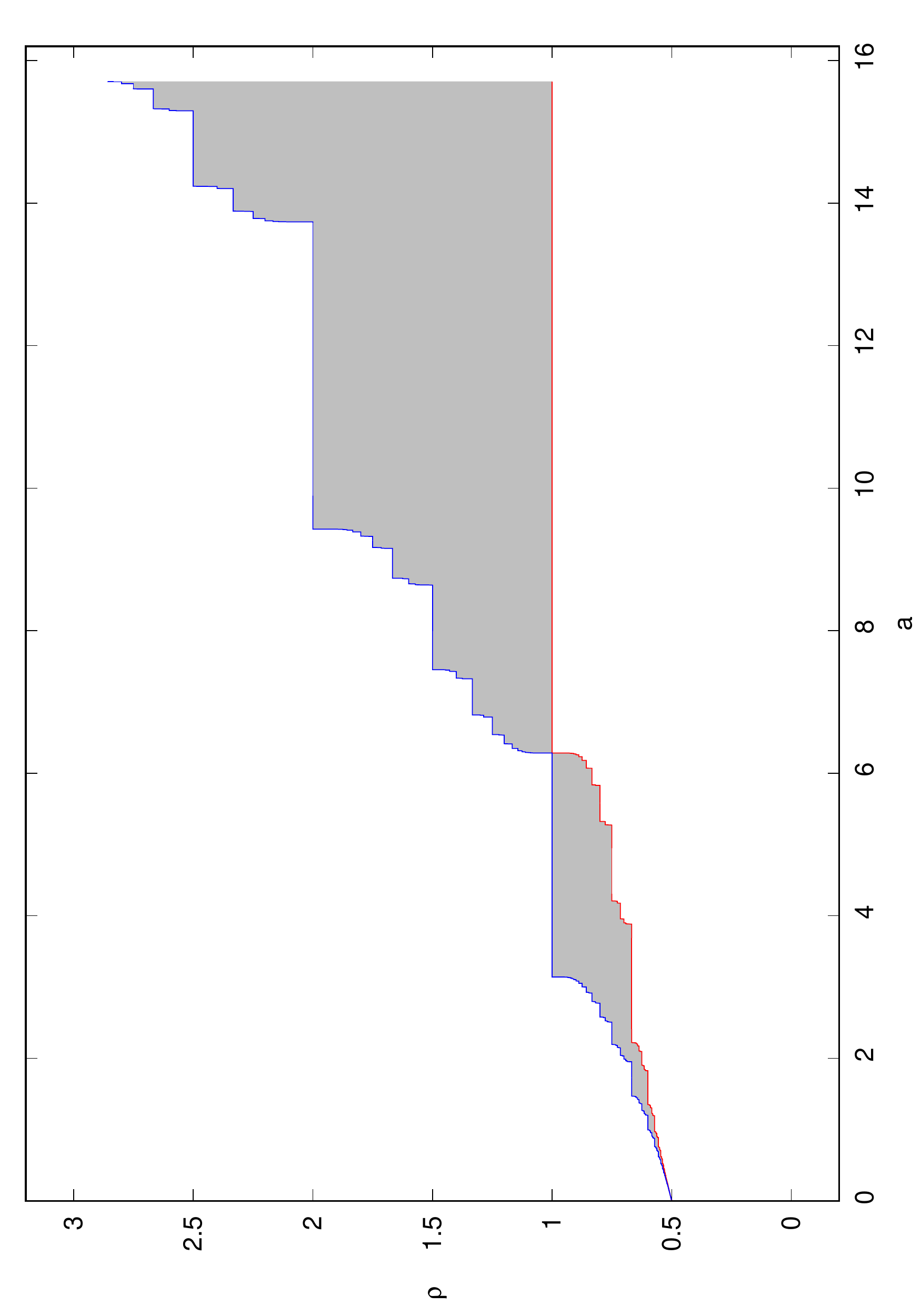}{Rotation interval graph for $\Omega = 0.5$.}
\hfill
\includeepssubfigure{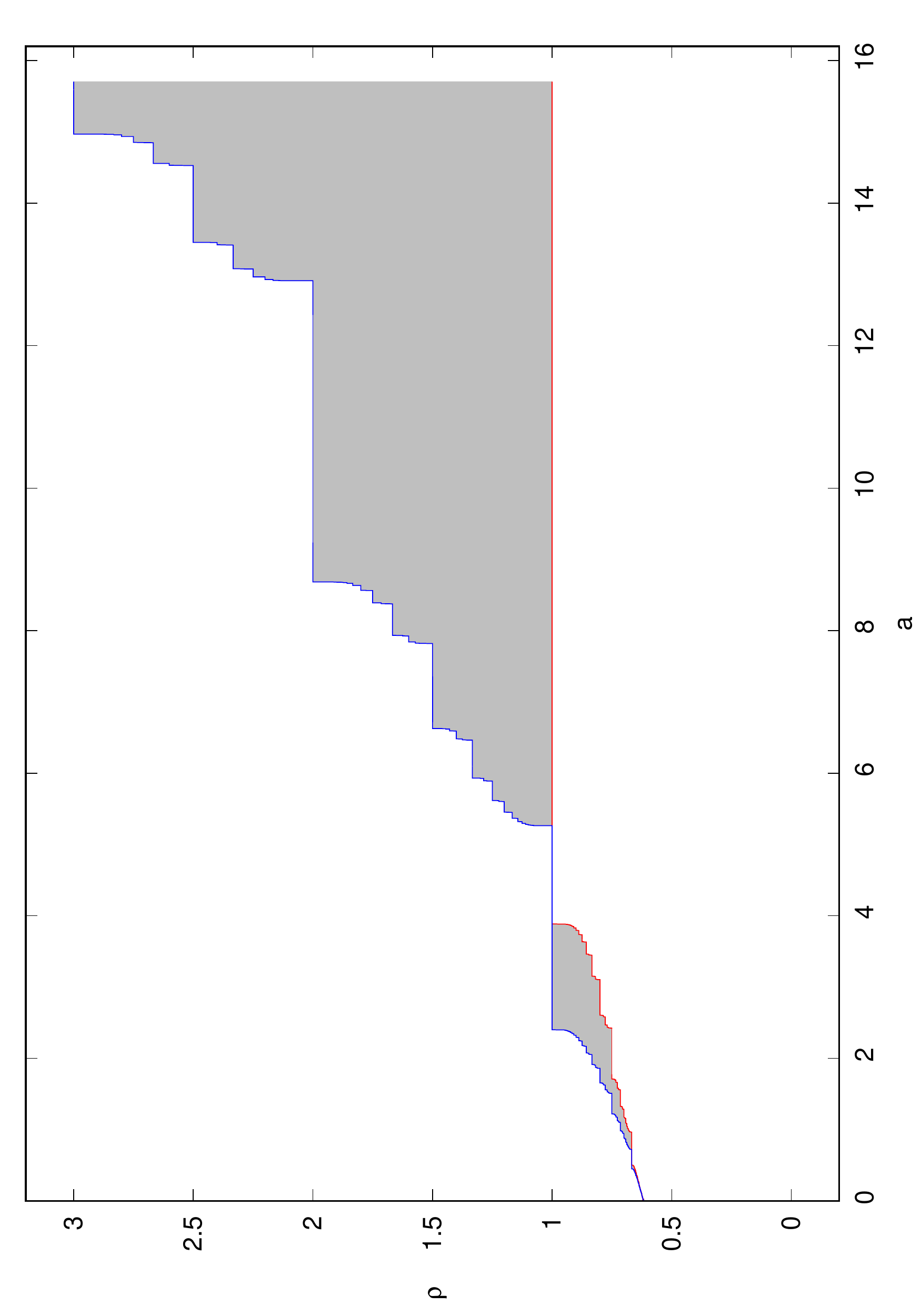}{Rotation interval graph with $\Omega$ equals to the Golden Mean.}
\par\par
\includepngsubfigure{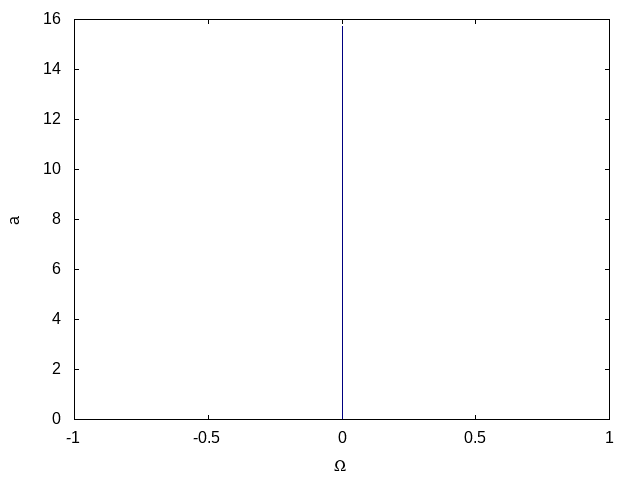}{$0-$Arnol'd ton\-gue.}
\hfill
\includepngsubfigure{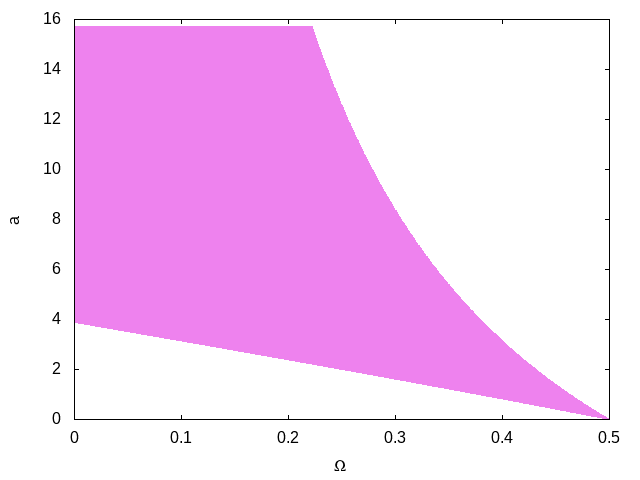}{$\mathcal{T}_{\varrho}$ Arnol'd tongue for $\varrho = 0.5$.}
\hfill
\includepngsubfigure{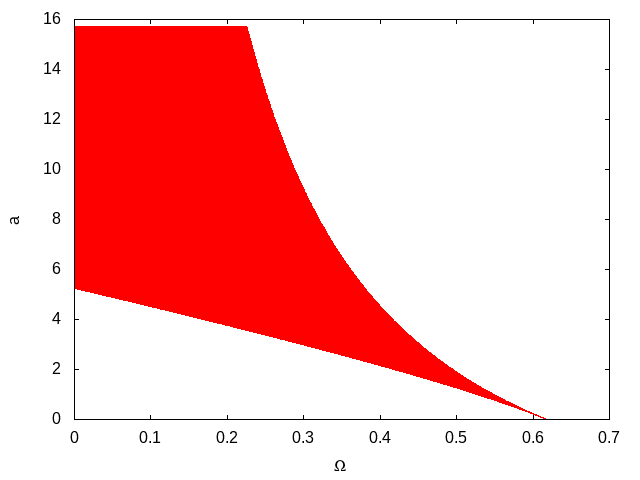}{$\mathcal{T}_{\varrho}$ for $\varrho$ equal to the Golden Mean.}
\captionsetup{width=\linewidth}
\caption{Graphs of the rotation interval and Arnol'd tongues for the discontinuous standard map
$\textsf{D}_{\Omega,a}$.
The graphs of the rotation intervals are plotted as a function of the parameter $a$.}\label{GraphsDiscStdMap}
\end{figure}

In Figure~\ref{GraphsDiscStdMap} we show some graphs of the rotation interval and
Arnol'd tongues for the discontinuous standard map.
The graphs of the rotation intervals are
plotted for three different values of $\Omega$ as a function of the parameter $a$.
The times taken for all the computation related with the rotation intervals and the Arnol'd Tongues for each of the families studied using Algorithms~\ref{alg:ClassicAlgo},~\ref{alg:SimosAlgo}~and~\ref{alg:CSBAlgo} can be found in Table~\ref{TimeTable}.
\section{Conclusions}\label{sec:Conclusions}

The algorithm proposed clearly outperforms all the other tested algorithms,
both in precision and speed even though the ``exact'' (and quick) part of the
algorithm does not work for all the non-decreasing liftings in $\dol$ which
have a constant section
(and hence the rotation number of these ``bad'' cases has to be computed with
the much more inefficient classical algorithm).
For all natural examples for which it has been tested, the
computational speed and precision were unparalleled. Moreover, the
set of functions becomes very general when one considers the fact
that the upper and lower functions inherently have constant sections
for any $F$ that is not strictly increasing. Hence, the
algorithm becomes a crucial tool to compute rotation intervals for
general functions in $\dol$ and hence to find the set of periods of
such maps \cite{ALM-Book}.\\
Moreover, a deeper study has been done on the dependence of the
rotation number on the parameters. Our preliminary results have found
that for irrational rotation numbers, the dependence of the
parameters around them is extremely sensitive, with continuity module
being at least $10^{25}.$ This agrees with Theorem~\ref{persistent},
which says that non-persistent functions have measure zero.

\bibliographystyle{plain}
\bibliography{NumRot-Algorithm}
\end{document}